\documentclass[11pt]{article} 
\usepackage{amsmath}
\usepackage{amsfonts}
\usepackage{amssymb}
\usepackage{stmaryrd}
\usepackage[utf8]{inputenc}
\usepackage[T1]{fontenc}
\usepackage[french]{babel}

\makeatletter
\def\@sect#1#2#3#4#5#6[#7]#8{%
  \ifnum #2>\c@secnumdepth
    \let\@svsec\@empty
  \else
    \refstepcounter{#1}%
    \protected@edef\@svsec{\@seccntformat{#1}\relax}%
  \fi
  \@tempskipa #5\relax
  \ifdim \@tempskipa>\z@
    \begingroup
      #6{%
        \@hangfrom{\hskip #3\relax\@svsec}%
          \interlinepenalty \@M #8\@@par}%
    \endgroup
    \csname #1mark\endcsname{#7}%
    \addcontentsline{toc}{#1}{%
      \ifnum #2>\c@secnumdepth \else
        \protect\numberline{\csname the#1\endcsname.}%
      \fi
      #7}%
  \else
    \def\@svsechd{%
      #6{\hskip #3\relax
      \@svsec #8}%
      \csname #1mark\endcsname{#7}%
      \addcontentsline{toc}{#1}{%
        \ifnum #2>\c@secnumdepth \else
          \protect\numberline{\csname the#1\endcsname.}%
        \fi
        #7}}%
  \fi
  \@xsect{#5}}
\def\@seccntformat#1{\csname the#1\endcsname.\quad}

\def\@begintheorem#1#2{\trivlist
   \item[\hskip \labelsep{\bfseries #1\ #2.}]\itshape}
\def\@opargbegintheorem#1#2#3{\trivlist
      \item[\hskip \labelsep{\bfseries #1\ #2\ (#3).}]\itshape}

\makeatother

\newtheorem{theo}[equation]{Th\'eor\`eme}
\newtheorem{lem}[equation]{Lemme}
\newtheorem{proposition}[equation]{Proposition}
\newtheorem{definition}[equation]{D\'efinition}
\newtheorem{cor}[equation]{Corollaire}
\newenvironment{remarque}{
\refstepcounter{equation}\trivlist%
\item[\hskip \labelsep{\bfseries Remarque \theequation.\ }]}%
{\endtrivlist}%

\makeatletter
\renewcommand\theequation{\thesection.\arabic{equation}}
\@addtoreset{equation}{section}
\makeatother
\newcommand{\carrenoir}{\rule{0.5em}{0.5em}}
\makeatletter
\newenvironment{demo}[1][\@empty]{\textbf{D\'emonstration~%
\ifx\@empty#1:\else #1~:\fi~}}
{\hfill\carrenoir\nolinebreak\vspace{2mm}}
\makeatother
\newcommand{\oper}[2]{\newcommand{#1}{\mathop{\mathrm{#2}}\nolimits} }
\oper{\Vol}{Vol}
\oper{\Hess}{Hess}
\oper{\Id}{Id}
\newcommand{\R}{\mathbb R}
\newcommand{\de}{\mathrm{ d }}
\oper{\Ker}{Ker}
\oper{\Ima}{Im}
\DeclareSymbolFont{greek}{OML}{ptmcm}{m}{it}
\DeclareMathSymbol{\codiff}{\mathord}{greek}{"0E}
\DeclareMathSymbol{\prodint}{\mathord}{greek}{"13}
\newcommand{\N}{\mathbb N}

\title{Sur la multiplicité des valeurs propres du laplacien
de Witten}
\author{Pierre Jammes}
\date{}
\begin{document}
\maketitle
{\small 
\textsc{Résumé.---}
Sur toute variété compacte de dimension supérieure ou égale à~4,
on prescrit le volume et le début du spectre du laplacien de Witten
agissant sur les $p$-formes différentielles pour $0<p<n$.
En particulier, on prescrit la multiplicité des premières valeurs propres.
Sur les variétés de dimension~3, on construit des exemples de première 
valeur propre multiple pour les 1-formes, dont la multiplicité dépend 
du genre maximal des surfaces immergées dont toute la 1-cohomologie est
induite par la cohomologie de la variété. En particulier, cette multiplicité
est au moins égale à~3.

Mots-clefs : laplacien de Witten, formes différentielles,
multiplicité de valeurs propres.

\medskip
\textsc{Abstract.---}
On any compact manifold of dimension greater than~4, we prescribe the volume
and any finite part of the spectrum of the Witten Laplacian acting on $p$-form
for $0<p<n$. In particular, we prescribe the multiplicity
of the first eigenvalues. On 3-dimensional manifolds, we give examples
of multiple first eigenvalue for 1-forms, whose multiplicity depends on
the maximal genus of embedded surfaces all of whose 1-cohomology is induced
by the cohomology of the manifold. In particular, this multiplicity is
at least~3.

Keywords : Witten Laplacian, differential forms, multiplicity of eigenvalues.

\medskip
MSC2000 : 58J50}

\section{Introduction}
Depuis que S.~Y.~Cheng a montré dans \cite{ch76} que la multiplicité de la
$k$-ième valeur propre du laplacien sur une surface compacte est majorée
en fonction de $k$ et de la topologie, ce problème de multiplicité 
pour le laplacien agissant sur les fonctions a fait l'objet de nombreux
travaux. En dimension~2, la majoration de Cheng, qui est aussi valable
pour les opérateur de Schrödinger, a été améliorée (voir~\cite{be80}, 
\cite{na88}, \cite{hohon99}), la meilleure estimation pour la multiplicité
de la 2\ieme~valeur propre d'un opérateur de Schrödinger ayant été obtenue 
par B.~Sévennec (\cite{se94}, 
\cite{se02}). On sait aussi que pour les opérateurs avec champ magnétique, 
la multiplicité des valeurs propres peut être arbitrairement grande
(\cite{cdvt93}, \cite{bcc98}, \cite{er02}). En dimension supérieure ou
égale à~3, Y. Colin de Verdière a montré (\cite{cdv86}, \cite{cdv87}) que 
toute rigidité disparaît et qu'on peut arbitrairement prescrire le début du
spectre avec multiplicité, et J.~Lohkamp a amélioré ce résultat en montrant
dans \cite{lo96} qu'on pouvait prescrire simultanément le
début du spectre, le volume et certains invariants de courbure.

Pour les opérateurs agissants sur les fibrés vectoriels naturels, 
ce problème est longtemps resté ouvert. Dans \cite{gu04}, P.~Guérini a 
montré qu'on peut prescrire toute partie finie du spectre du laplacien 
de Hodge-de~Rham, qui agit sur les formes différentielles, mais en 
imposant aux valeurs propres prescrites d'être simples. Un résultat
semblable a été obtenu par M.~Dahl pour l'opérateur de Dirac (\cite{da05}).

Un premier résultat de multiplicité est obtenu dans \cite{ja09},
où on construit un nombre arbitraire de valeurs propres doubles pour
le laplacien de Hodge-de~Rham. Plus récemment, j'ai montré dans
\cite{ja11} qu'on pouvait étendre le résultat de Colin de Verdière
en prescrivant le début du spectre du laplacien de Hodge-de~Rham avec
multiplicité sur les variétés de dimension supérieure ou égale à~6,
mais la technique échoue dans le cas des petites dimensions, des formes de 
degré $[\frac n2]$ et de la première valeur propres des 1-formes.

Le but de cet article est de montrer qu'en dimension supérieure ou
égale à~4, on peut prescrire le spectre du laplacien de Witten sans
restriction sur le degré. Cet opérateur, agissant sur les formes
différentielles, a été popularisé par E.~Witten dans \cite{wi82},
où il l'utilise entre autres pour redémontrer les inégalités de Morse
(voir par exemple~\cite{he84}). Contrairement à ce qui se fait
habituellement, nous n'étudierons pas ici la limite semi-classique de cet
opérateur. On utilisera plutôt le fait que l'étude du spectre
du laplacien de Witten revient à découpler la métrique et la mesure
(cette notion sera précisée dans la section~\ref{theo}).
Cet opérateur peut être vu comme un analogue de l'opérateur de Schrödinger
pour les formes différentielles: dans le cas des formes de degré~0,
c'est-à-dire des fonctions, on retrouve tous les opérateurs de Schrödinger
dont la première valeur propre est nulle.

Rappelons la définition de cet opérateur (voir la section~\ref{theo}
pour plus de détails). Étant données une variété compacte
$M^n$ et  une fonction $\varphi\in C^\infty$, on définit une différentielle
tordue par $\tilde\de\omega=\de\omega+\de\varphi\wedge\omega$ et on
note $\tilde\codiff$ son adjoint. Le laplacien de Witten est alors
défini par $\tilde\Delta_\varphi=\tilde\de\tilde\codiff+
\tilde\codiff\tilde\de$. On retrouve la laplacien de Hodge quand $\varphi$
est constante. La théorie de Hodge s'applique au laplacien de Witten,
en particulier on a une décomposition de Hodge $\Omega^p(M)=\tilde\de
\Omega^{p-1}(M)\oplus\Ker\tilde\Delta_\varphi\oplus
\tilde\codiff\Omega^{p+1}(M)$, qui est stable par le laplacien. Si
on note
\begin{equation}
0<\tilde\mu_{p,1}(M,g,\varphi)\leq\tilde\mu_{p,2}(M,g,\varphi)\leq\ldots
\end{equation}
Les valeurs propres du laplacien de Witten en restriction à 
$\tilde\codiff\Omega^{p+1}(M)$, alors son spectre en restriction
à $\tilde\de \Omega^{p-1}(M)$ est $(\tilde\mu_{p-1,i}(M,g,\varphi))$. 
Le spectre non nul du laplacien de Witten se déduit donc des 
$\tilde\mu_{p-1,i}(M,g,\varphi)$, $0\leq p\leq n-1$.

Avec ces notations, le premier résultat de cet article peut s'énoncer 
comme suit: 
\begin{theo}\label{intro:th1}
Soit $M^n$ une variété compacte connexe orientable (avec ou sans bord) 
de dimension 
$n\geq4$ et $N\in\N^*$. Si on se donne un réel $V>0$ et et des suites 
$0<a_{p,1}\leq a_{p,2}\leq\ldots\leq a_{p,N}$
pour $1\leq p\leq n-1$, alors il existe une métrique $g$ et une fonction
$\varphi$ sur $M$ telles que
\begin{itemize}
\item $\tilde\mu_{p,k}(M,g,\varphi)=a_{p,k}$ pour $1\leq k\leq N$ et
$1\leq p\leq n-1$;
\item $\Vol(M,g)=V$.
\end{itemize}
\end{theo}
On peut donc en particulier prescrire la multiplicité des premières
valeurs propres pour tous les degrés non nuls simultanément.

Comme dans \cite{ja11}, certaines techniques échouent en dimension~3.
On peut cependant construire une première valeur propre dont la
multiplicité dépend de la topologie, étendant ainsi un résultat 
que Colin de Verdière avait obtenu sur les surfaces. Étant donnée
une surface compacte $\Sigma$, et en notant $C(\Sigma)$ le nombre
chromatique de $\Sigma$, c'est-à-dire le plus grand entier $k$ tel
que le graphe complet à $k$ sommets se plonge dans $\Sigma$, il montre
dans \cite{cdv87} qu'il existe un opérateur de Schrödinger sur 
$\Sigma$ dont la multiplicité de la 2\ieme{} valeur propre est
$C(\Sigma)-1$. Le résultat qu'on peut obtenir en dimension~3 est
le suivant:
\begin{theo}\label{intro:th2}
Soit $(M,g)$ une variété riemannienne compacte de dimension~3. Pour
toute surface plongée $\Sigma\hookrightarrow M$ telle que l'application
naturelle $H^1(M)\to H^1(\Sigma)$ en cohomologie de~Rham soit surjective,
il existe une métrique $g$ et une fonction $\varphi$ sur $M$ telles que 
$\tilde\mu_{1,1}(M,g,\varphi)$ soit de multiplicité $C(\Sigma)-1$.
\end{theo}
\begin{remarque}
La condition de surjectivité de $H^1(M)\to H^1(\Sigma)$ peut s'interpréter
comme un analogue cohomologique de la notion de surface incompressible.
\end{remarque}
Le seul résultat de multiplicité en dimension~3 pour le laplacien de
Hodge-de~Rham était jusqu'à présent la construction de valeurs
propres doubles donnée dans \cite{ja09}. En appliquant le 
théorème~\ref{intro:th2} avec $\Sigma=S^2$, on obtient qu'on
peut toujours faire mieux pour la première valeur propre du laplacien
de Witten :
\begin{cor}
Soit $(M,g)$ une variété riemannienne compacte de dimension~3.
Il existe une métrique  $g$ et une fonction $\varphi$ sur $M$ telles que
$\tilde\mu_{1,1}(M,g,\varphi)$ soit de multiplicité~3.
\end{cor}

 L'article est organisé comme suit : la section~\ref{theo} est consacrée
à des rappels sur la théorie du laplacien de Witten et à la démonstration
de quelques lemmes techniques. On verra en particulier qu'on peut
définir son spectre sans introduire la différentielle tordue en remplaçant
la mesure riemannienne $\de v_g$ par la mesure $e^{-2\varphi}\de v_g$
(lemme~\ref{theo:spectre3}). On montrera dans la section~\ref{minoration}
un résultat de rigidité conforme du spectre qui généralise celui 
obtenu dans \cite{ja07b} pour le laplacien de Hodge et qui sera utilisé
dans la section suivante. Dans les sections~\ref{conv} et~\ref{boule} 
on donnera deux résultats de convergence spectrale pour le laplacien
de Witten. Le premier traite de la convergence du spectre d'une variété
vers celui d'un de ses domaines, généralisant des résultats
obtenus pour le laplacien agissant sur les fonctions \cite{cdv86} et 
le laplacien de Hodge \cite{ja11}.
Le second concerne les variétés privées d'une petite boule et généralise
un théorème de C.~Anné et B.~Colbois \cite{ac93}. Enfin, on démontrera
dans la dernière section les théorème~\ref{intro:th1} et~\ref{intro:th2} en
utilisant les deux théorèmes de convergence spectrale.

Ce travail a été mené avec le soutien du projet ANR Geodycos.

\section{Théorie de Hodge et conditions de bord pour le laplacien de Witten}%
\label{theo}
\subsection{Laplacien de Witten}
On va rappeler dans ce paragraphe la définition et les propriétés 
élémentaires du laplacien de Witten. On peut se référer à \cite{wi82} ou
\cite{hn06} pour une présentation plus complète.

Étant donnée une fonction $\varphi\in C^\infty(M)$, on définit une 
différentielle tordue $\tilde\de:\Omega^p(M)\to\Omega^{p+1}(M)$ par
\begin{equation}\label{theo:eq1}
\tilde\de\omega=e^{-\varphi}\de(e^{\varphi}\omega)=
\de\omega+\de\varphi\wedge\omega
\end{equation}
et une codifférentielle tordue $\tilde\codiff:\Omega^p(M)\to\Omega^{p-1}(M)$ 
par
\begin{equation}\label{theo:eq2}
\tilde\codiff\omega=e^{\varphi}\codiff(e^{-\varphi}\omega)=
\codiff\omega+\prodint_{\nabla\varphi}\omega.
\end{equation}
On peut vérifier que ces deux opérateurs sont adjoints l'un de l'autre.
Le laplacien de Witten associé à $\varphi$ est alors défini par
\begin{equation}\label{theo:eq3}
\tilde\Delta_\varphi=\tilde\de\tilde\codiff+\tilde\codiff\tilde\de
:\Omega^p(M)\to\Omega^p(M).
\end{equation}

On peut réécrire le laplacien de Witten sous les formes suivantes :
\begin{theo}\label{theo:th1}
Pour toute fonction $\varphi$, on a
\begin{eqnarray*}
\tilde\Delta_\varphi&=&\Delta+|\de\varphi|^2+\mathcal L_{\nabla\varphi}+
\mathcal L_{\nabla\varphi}^*\\
&=&\Delta+|\de\varphi|^2+\Delta\varphi-(\mathcal L_{\nabla\varphi}g_p)\\
&=&\Delta+|\de\varphi|^2+\Delta\varphi+2(\Hess\varphi).
\end{eqnarray*}
\end{theo}
Précisons les notations de cet énoncé : $\mathcal L_{\nabla\varphi}$ 
désigne la dérivée de Lie
par rapport à $\nabla\varphi$, $\mathcal L_{\nabla\varphi}^*$ son adjoint 
$L^2$ et $g_p$ la métrique induite par $g$ sur le fibré des $p$-formes.
Par conséquent, $\mathcal L_{\nabla\varphi}g_p$ est ponctuellement
une forme quadratique sur $\Lambda^pTM$ qu'on identifie à l'endomorphisme
symétrique de $\Lambda^pTM$ canoniquement associé. Le hessien 
$\Hess\varphi$ de $\varphi$ est une forme bilinéaire symétrique sur les
champs de vecteur, donc sur les $1$-formes, qu'on identifie là encore
à un endomorphisme de $\Lambda^1TM$. On l'étend à $\Lambda^pTM$ de la 
manière suivante : Si $(e_1,\ldots,e_n)$ est une base orthonormée de 
$\Lambda^1TM$, alors
\begin{equation}\label{theo:eq4}
(\Hess\varphi)(e_1\wedge\ldots\wedge e_p)=\sum_ie_1\wedge\ldots\wedge
e_{i-1}\wedge\Hess\varphi(e_i)\wedge e_{i+1}\wedge\ldots\wedge e_p.
\end{equation}

La première égalité du théorème~\ref{theo:th1} s'obtient en développant
l'expression (\ref{theo:eq3}). La démonstration des autres équations est 
rarement détaillée dans la littérature, nous la rappelons en appendice.

En outre, le laplacien de Witten vérifie la relation de commutation suivante
avec le dualité de Hodge:
\begin{equation}\label{theo:hodge}
\tilde\Delta_\varphi*=*\tilde\Delta_{-\varphi}
\end{equation}
Cette relation implique que $\tilde\mu_{p,i}(g,\varphi)=
\tilde\mu_{n-p-1,i}(g,-\varphi)$.

Le spectre du laplacien de Hodge-de~Rham sur le produit riemannien de
deux variétés peut se calculer à l'aide de la formule de Künneth. On 
démontre ci-dessous la généralisation de cette formule au laplacien 
de Witten, qui nous sera utile pour la construction de valeurs propres
multiples.
\begin{theo}[Formule de Künneth]\label{theo:kunneth}
Soit $(M_1,g_1)$ et $(M_2,g_2)$ deux variétés riemanniennes, 
$\varphi_i\in C^\infty(M_i)$ deux fonctions et $\alpha_i\in\Omega(M_i)$ 
deux formes différentielles sur les variétés $M_1$ et $M_2$. Alors on 
a sur le produit riemannien $(M_1\times M_2,g_1\oplus g_2)$
$$\tilde\Delta_\varphi(\alpha_1\wedge\alpha_2)=\tilde\Delta_{\varphi_1}
\alpha_1\wedge\alpha_2+\alpha_1\wedge\tilde\Delta_{\varphi_2}\alpha_2$$
où $\varphi\in C^\infty(M_1\times M_2)$ est définie par $\varphi=
\varphi_1+\varphi_2$, les formes $\alpha_i$ et les fonctions $\varphi_i$ 
étant identifiées à leur relevé sur $M_1\times M_2$.
\end{theo}
\begin{demo}
On va exploiter l'une des expressions du laplacien de Witten données par 
le théorème~\ref{theo:th1}.

On sait que pour le laplacien de Hodge-de~Rham, on a
\begin{equation}\label{theo:ku1}
\Delta(\alpha_1\wedge\alpha_2)=\Delta\alpha_1\wedge\alpha_2+\alpha_1
\wedge\Delta\alpha_2.
\end{equation}
Comme $\de\varphi_1$ et $\de\varphi_2$ sont orthogonaux pour 
la métrique $g_1\oplus g_2$, on a aussi 
\begin{equation}\label{theo:ku2}
|\de\varphi|^2=|\de\varphi_1|^2+|\de\varphi_2|^2.
\end{equation}
 Enfin, chaque fonction $\varphi_i$ ne varie que dans la direction de $M_i$,
donc le hessien $\Hess\varphi$ agissant sur $\Omega^1(M_1\times M_2)$
prend la forme
\begin{eqnarray}
\Hess\varphi & = & \Hess\varphi_1\oplus\Id_{\Lambda^1TM_2}
+\Id_{\Lambda^1TM_1}\oplus \Hess\varphi_2.
\end{eqnarray}
Étendu au $p$-formes à l'aide de (\ref{theo:eq4}), s'écrit
\begin{equation}\label{theo:ku4}
(\Hess\varphi)(\alpha_1\wedge\alpha_2)=
(\Hess\varphi_1)\alpha_1\wedge\alpha_2+
\alpha_1\wedge(\Hess\varphi_2)\alpha_2
\end{equation}

 On obtient la formule souhaitée en sommant (\ref{theo:ku1}), 
(\ref{theo:ku2}), et (\ref{theo:ku4}).
\end{demo}

\subsection{Cohomologie}
La différentielle tordue vérifie $\tilde\de^2=0$ et la théorie de Hodge
s'applique à $\tilde\de$. Le noyau de $\tilde\Delta_\varphi$ en restriction 
aux $p$-formes est donc isomorphe à la cohomologie de $\tilde\de$. Il
s'avère que la dimension de cette cohomologie ne dépend pas de $\varphi$,
cela découle de la remarque suivante :
\begin{lem}\label{theo:T}
Si on pose $T_\varphi\omega=e^\varphi\omega$ alors $T_\varphi\tilde\de=
\de T_\varphi$. En particulier, l'application $T_\varphi$ est un
isomorphisme du noyau (resp. l'image) de $\tilde\de$ vers le noyau
(resp. l'image) de $\de$.
\end{lem}
Ce lemme interviendra plusieurs fois par la suite, et on en déduit
dès à présent :
\begin{cor}\label{theo:coh}
 Les cohomologies de $\tilde\de$ et $\de$ sont isomorphes. Si $h$ est
une $p$-forme $\tilde\Delta_\varphi$-harmonique, alors $h$ minimise
la norme $L^2$ pour la mesure $\de v_g$ dans sa classe de 
$\tilde\de$-cohomologie. De plus, $T_\varphi h$ est $\de$-fermée et 
minimise la norme $L^2$ pour la mesure $\de v_\varphi=e^{-2\varphi}\de v_g$ 
dans sa classe de $\de$-cohomologie.
\end{cor}
\begin{demo}
 L'isomorphie entre les cohomologies découle immédiatement du 
lemme~\ref{theo:T}.

Le fait que $h$ minimise la norme $L^2$ pour la mesure $\de v_g$ résulte
du fait que $\tilde\codiff h=0$ et donc que $h$ est orthogonale au
forme $\tilde\de$-exacte.

Enfin, $T_\varphi h$ est fermée d'après le lemme~\ref{theo:T}, et si
$\alpha$ est une $(p-1)$-forme, alors 
$\int_M\langle T_\varphi h,\de\alpha\rangle \de v_\varphi=
\int \langle e^{-\varphi}T_\varphi h,e^{-\varphi}\de\alpha\rangle\de v_g
=\int\langle h,\tilde\de(e^{-\varphi}\alpha)\rangle\de v_g=0$.
Par conséquent, $T_\varphi h$ est orthogonale pour la mesure $\de v_\varphi$
au formes $\de$-exacte, donc minimise la norme $L^2$ pour la mesure
$\de v_\varphi$ dans sa classe de $\de$-cohomologie.
\end{demo}

\subsection{Variétés à bord}
Les énoncés des théorèmes de convergence spectrale et les démonstrations 
des résultats de multiplicité font intervenir des 
variétés à bord. Nous rappelons ou démontrons ici quelques propriétés
du laplacien de Witten dans ce cadre.

Si $U$ est un domaine à bord $C^1$ d'une variété compacte $M$,
on note $j:\partial U\to \overline U$
l'injection canonique et $N$ un champ de vecteur normal au bord.
Les conditions de bord classiques du laplacien de Hodge-de~Rham se
généralisent au laplacien de Witten. La condition absolue 
($\mathrm A_\varphi$) s'écrit
\begin{equation}
(\mathrm A_\varphi)\ \left\{\begin{array}{l}j^*(\prodint_N\omega)=0\\
j^*(\prodint_N\tilde\de\omega)=0\end{array}\right.\textrm{ ou }
\left\{\begin{array}{l}j^*(*\omega)=0\\
j^*(*\tilde\de\omega)=0\end{array}\right.
\end{equation}
et la condition relative ($\mathrm R_\varphi$) est définie par
\begin{equation}
(\mathrm R_\varphi)\ \left\{\begin{array}{l}j^*(\omega)=0\\
j^*(\tilde\codiff\omega)=0\end{array}\right.
\end{equation}
Quand $\varphi=0$, on retrouve les conditions absolues et relatives
usuelles du laplacien de Hodge.

La dualité de Hodge transforme la condition ($\mathrm A_\varphi$) en la 
condition ($\mathrm R_{-\varphi}$), ce qui est cohérent avec
la formule~\ref{theo:hodge}. B.~Heffler et F.~Nier montrent dans \cite{hn06}
que la condition ($\mathrm R_\varphi$) est admissible, c'est-à-dire que 
l'opérateur associé est elliptique, il en va donc de même pour la condition 
($\mathrm A_\varphi$). En outre, L'opérateur $T:\omega\to e^\varphi\omega$ 
défini dans le lemme~\ref{theo:T} transforme la condition 
($\mathrm A_\varphi$) en la condition ($\mathrm A_0$)

Le corollaire~\ref{theo:coh} est valable dans ce contexte. En particulier,
si $h$ est une forme $\tilde\Delta_\varphi$ harmonique vérifiant la
condition ($\mathrm A_\varphi$), alors $T_\varphi h$ est fermée, vérifie
la condition ($\mathrm A_0$) et minimise la norme $L^2$ pour la mesure
$\de v_\varphi$ dans sa classe de cohomologie.

On utilisera aussi la condition suivante, qui généralise la condition
de Dirichlet pour les fonctions :
\begin{equation}
(D)\ \left\{\begin{array}{l}j^*(\omega)=0\\
j^*(*\omega)=0\end{array}\right.
\end{equation}
Elle n'interviendra que pour le laplacien de Hodge, et
on utilisera principalement le fait que pour cette condition, le noyau
du laplacien est trivial (voir~\cite{an89}), et donc que son spectre
est strictement positif.

Pour finir ce paragraphe nous allons donner la généralisation suivante 
d'un résultat connu pour la différentielle standard. Cette propriété 
interviendra en particulier au paragraphe suivant pour donner une
caractérisation variationnelle du spectre sur les variétés à bord 
(propositions~\ref{theo:spectre1} à \ref{theo:spectre3}) :

\begin{theo}\label{theo:prim}
Si $\omega$ est une forme exacte pour la différentielle tordue $\tilde\de$
sur $\overline U$, alors il existe une forme $\psi$ vérifiant la condition 
($A_\varphi$) et 
telle que $\omega=\tilde\de\tilde\codiff\psi$. En particulier, la forme 
$\theta=\tilde\codiff\psi$ vérifie $j^*(\prodint_N\theta)=0$, est 
orthogonale aux formes $\alpha$ telles que $\tilde\de\alpha=0$ et minimise 
la norme $L^2$ dans $\tilde\de^{-1}\omega$.

En outre, la forme $\omega$ est contenue dans l'adhérence $L^2$ des formes 
exactes vérifiant la condition $(A_\varphi)$.
\end{theo}
\begin{remarque}\label{theo:rem}
Le fait que $\theta$  minimise la norme $L^2$ dans $\tilde\de^{-1}\omega$
est équivalent au fait que $T_\varphi\theta$, qui vérifie la condition
de bord ($A_0$), minimise la norme $L^2$ relative à la mesure 
$e^{-2\varphi}\de v_g$ dans $\de^{-1}T_\varphi\omega$.
\end{remarque}

Grâce au lemme qui suit, la démonstration est exactement la même
que pour la différentielle standard (voir~\cite{ja11}, proposition~2.4).
\begin{lem}[intégration par partie]\label{theo:ipp}
Si $\alpha\in\Omega^p(U)$ et $\beta\in\Omega^{p+1}(U)$ alors
$$(\tilde\de\alpha,\beta)=(\alpha,\tilde\codiff\beta)+\int_{\partial U}
j^*(\alpha\wedge*\beta).$$
\end{lem}
\begin{demo}
D'une part, on a la formule classique $(\de\alpha,\beta)=
(\alpha,\codiff\beta)+\int_{\partial U}j^*(\alpha\wedge*\beta)$ qui découle
de la formule de Stokes. D'autre part, le produit intérieur par un vecteur
est l'adjoint du produit extérieur par la 1-forme duale,
c'est-à-dire qu'en tout point $x\in U$, $\langle\de\varphi\wedge\alpha,
\beta\rangle_x=\langle\alpha,\prodint_{\nabla\varphi}\beta\rangle_x$.
Comme cette formule est ponctuelle, son intégrale sur $U$ ne
fait pas apparaître de terme de bord. L'addition des deux équations
donne le résultat souhaité.
\end{demo}
\subsection{Caractérisations du spectre du laplacien de Witten}
L'étude du spectre du laplacien de Hodge-de~Rham est grandement
facilitée par le principe variationnel suivant qui remonte à J.~Dodziuk
et qui reste valable pour le laplacien de Witten :
\begin{proposition}[\cite{do82}, \cite{mc93}]\label{theo:spectre1}
 Sur une variété compacte sans bord ou avec condition de bord $(A_\varphi)$, 
on a
$$\tilde\mu_{p,i}=\inf_{V_i}\sup_{\omega\in V_i\backslash\{0\}}\left\{
\frac{\|\omega\|^2}{\|\theta\|^2},\ \tilde\de\theta=\omega\right\},$$
où $V_i$ parcourt l'ensemble des sous-espaces de dimension $i$ dans
l'espace des $p+1$-formes exactes lisses.
\end{proposition}
 Ici encore, comme la théorie de Hodge fonctionne pour $\tilde\de$, la
démonstration est la même que dans le cas du laplacien de Hodge-de~Rham.
 Dans le cas à bord, le fait qu'on obtienne le 
spectre pour la condition $(A_\varphi)$ même si on exige aucune condition
sur $\omega$ et $\theta$ découle du théorème~\ref{theo:prim} 
(en particulier du fait que l'adhérence des formes exactes vérifiant 
$(A_\varphi)$ contient toutes les formes exactes de $\overline U$).

Comme remarqué dans \cite{ja11}, on peut 
reformuler ce résultat de la manière suivante:
\begin{proposition}\label{theo:spectre2}
Le spectre et les espaces propres du laplacien de Witten en restriction à
$\Ima\tilde\de$ sont ceux de la forme quadratique 
$Q(\omega)=\|\omega\|^2_{L^2}$ relativement à la norme 
$|\omega|=\displaystyle\inf_{\tilde\de\theta=\omega}
\|\theta\|_{L^2}$.
\end{proposition}
Enfin, grâce au lemme~\ref{theo:T}, on peut donner cette autre formulation 
qui fait appel à la différentielle $\de$ au lieu de $\tilde\de$, et qui
montre qu'on peut interpréter le laplacien de Witten par un découplage entre
la métrique et la mesure.
C'est sous cette forme qu'on utilisera cette caractérisation du spectre.
\begin{proposition}\label{theo:spectre3}
Le spectre et les espaces propres, transportés par $T_\varphi$, du laplacien 
de Witten en restriction à $\Ima\tilde\de$ sont ceux de la forme 
quadratique $Q(\omega)=\int_M |\omega|^2\de v_\varphi$ définie sur les
formes exactes, relativement à la norme
$|\omega|=\displaystyle\inf_{\de\theta=\omega}\int_M|\theta|^2\de v_\varphi$,
où $\de v_\varphi$ désigne la mesure $e^{-2\varphi}\de v_g$.
\end{proposition}
\begin{demo}
Il suffit de remarquer que si on a $\tilde\de\tilde\theta=\tilde\omega$, 
et qu'on pose $\omega=T_\varphi\tilde\omega$ et 
$\theta=T_\varphi\tilde\theta$, alors
$T_\varphi\tilde\de\tilde\theta=\de T_\varphi\tilde\theta=
T_\varphi\tilde\omega$, donc $\de\theta=\omega$. Par ailleurs, on a aussi
$\int |\tilde\omega|^2\de v_g=\int |\omega|^2e^{-2\varphi}\de v_g$ et
$\int |\tilde\omega|^2\de v_g=\int |\omega|^2e^{-2\varphi}\de v_g$.
\end{demo}

\begin{remarque}
Comme dans le cas du laplacien de Hodge-de~Rham, on peut déduire de
la proposition~\ref{theo:spectre3} que le spectre du laplacien de Witten
est continu pour la topologie $C^0$ sur l'espace des métriques. On peut
aussi en déduire la continuité $C^0$ du spectre par rapport à $\varphi$,
ce qui n'était pas du tout évident dans les expressions du 
théorème~\ref{theo:th1}.
\end{remarque}

\section{Minoration du spectre dans une classe conforme pondérée
et inégalités de Sobolev}\label{minoration}

 Le but de cette section est de généraliser au laplacien de Witten
une minoration conforme du spectre du laplacien de Hodge-de~Rham
obtenue dans \cite{ja07b}. Plus précisément, étant donné une
métrique $g$, une fonction $\varphi$ et un réel $\alpha>0$, on définit
la classe conforme pondérée de poids $\alpha$ de $(g,\varphi)$ par
\begin{equation}
[g,\varphi]_\alpha=\left\{(e^{2u}g,\varphi-\alpha u),\ 
u\in C^\infty(M)\right\}.
\end{equation}
On va montrer que pour certains degrés $p$, dépendants de $\alpha$, la valeur
propre $\tilde\mu_{p,1}(M,g,\varphi)$ est uniformément minorée sur 
$[g,\varphi]_\alpha$, à volume fixé. Le cas traité dans \cite{ja07b} 
correspond à $\varphi=0$, $\alpha=0$ et $p\in[\frac n2-1,\frac n2]$.
\begin{theo}\label{min:theo}
Si $(M^n,g)$ est une variété riemannienne compacte de dimension~$n$, 
$\varphi$ une fonction sur $M$, $\alpha>0$ un réel et $p\in[1,n]$ un entier.
Si $p\in[\frac n2-\alpha-1,\frac n2-\alpha]$. Alors il existe 
une constante $c>0$ ne dépendant que de la classe conforme pondérée 
$[g,\varphi]_\alpha$ telle que $$\tilde\mu_{p,1}(M,g,\varphi)
\Vol(M,g)^{\frac 2n}\geq c.$$
\end{theo}

Ce théorème repose sur l'existence d'inégalités de Sobolev pour les
formes différentielles. Ces inégalités sont apparues sous la plume de 
V.~Gol'dshtein et M.~Troyanov \cite{gt06}, et certaines des constantes
qu'elles font intervenir
sont des invariants conformes. Nous allons ici en donner
une version «~à poids~» qui permettra l'application au laplacien de
Witten. On pourrait déduire ces inégalités des résultats de \cite{gt06},
mais par soucis d'exhaustivité ou va les montrer d'une manière plus directe 
en utilisant une méthode qui était déjà esquissée dans \cite{ja07a}. Avant 
de démontrer le théorème~\ref{min:theo}, nous allons donc donner deux 
lemmes, l'un sur l'existence d'inégalités de Sobolev et l'autre
sur l'invariance conforme pondérée de certaines des constantes de Sobolev 
associées. Par commodité, pour un réel $r>1$, on notera 
$\|\omega\|_{r,\varphi}$ la norme $L^p$ de $\omega$ pour la mesure 
$e^{-r\varphi}\de v_g$.

\begin{lem}\label{min:lem1}
Soit $(M^n,g)$ est une variété riemannienne compacte de dimension~$n$, 
$\varphi$ une fonction sur $M$ et $p\in[1,n]$ un entier. Si $r,s>1$ sont
deux réels tels que $\frac1s-\frac1r\leq\frac1n$ il existe 
une constante $K>0$ dépendant de $g$, $p$, $\varphi$ $r$ et $s$ telle que
$$ \inf_{\de\theta=0}\|\omega-\theta\|_{r,\varphi}\leq 
K\|\de\omega\|_{s,\varphi}.$$
\end{lem}
\begin{demo}
L'opérateur $\tilde\de+\tilde\codiff$ est elliptique, donc il existe des
constantes $A>0$ et $B>0$ dépendant de $g$ et $\varphi$ telles que si 
$\omega\in\Omega^p(M)$ est orthogonale à son noyau, alors 
$\|\nabla\omega\|_s\leq A\|\tilde\de\omega\|_s+B\|\tilde\codiff\omega\|_s$. 
De plus, toujours en supposant que $\omega$ est orthogonale au noyau de 
$\tilde\de+\tilde\codiff$ (donc aux formes parallèles), on a l'inégalité
de Sobolev $\|\omega\|_r\leq A'\|\nabla\omega\|_s$ pour une constante
$A'>0$ ne dépendant que de $g$. Par 
théorie de Hodge, il existe une forme $\theta\in\Omega^p(M)$ telle que 
$\tilde\de\theta=0$ et $\tilde\codiff(\omega-\theta)=0$. En appliquant les 
inégalités précédentes à $\omega-\theta$ on obtient
\begin{equation}
\|\omega-\theta\|_r\leq A'A\|\tilde\de\omega\|_s.
\end{equation}
Si on pose $\omega'=T_\varphi\omega$ et $\theta'=T_\varphi\theta$,
la condition $\tilde\de\theta=0$ devient $\de\theta'=0$ et les normes
précédentes s'écrivent
\begin{equation}
\|\omega-\theta\|_r=
\left(\int_Me^{-r\varphi}|\omega'-\theta'|^r\de v_g\right)^{\frac1r}
=\|\omega'-\theta'\|_{r,\varphi}
\end{equation}
et
\begin{equation}
\|\tilde\de\omega\|_r=
\left(\int_Me^{-s\varphi}|\de\omega'|^s\de v_g\right)^{\frac1s}
=\|\de\omega'\|_{s,\varphi}.
\end{equation}
On a donc montré l'existence d'une constante $K=A'A>0$ telle que pour toute
forme $\omega'$, il existe une forme fermée $\theta'$ telle que 
$\|\omega'-\theta'\|_{r,\varphi}\leq K\|\de\omega'\|_{s,\varphi}$.
L'inégalité du lemme s'en déduit immédiatement.
\end{demo}

\begin{lem}\label{min:lem2}
Si $(M^n,g)$ est une variété riemannienne compacte de dimension~$n$, 
$\varphi$ une fonction sur $M$ et $p\in[1,n]$ un entier, et $\alpha$ un réel. 
alors la constante
$$\sup_{\omega\in\Omega^p(M)}
\inf_{\de\theta=0}\frac{\|\omega-\theta\|_{\frac n{p-\alpha},\varphi}}%
{\|\de\omega\|_{\frac n{p+1-\alpha},\varphi}}$$
est strictement positive et est un invariant de la classe conforme pondérée 
$[g,\varphi]_\alpha$ qu'on notera $K_p(M,[g,\varphi]_\alpha)$.
\end{lem}
\begin{demo}
Si on pose $r=\frac n{p-\alpha}$ et $s=\frac n{p+1-\alpha}$, on a
$\frac1s-\frac1r=\frac1n$ et on peut donc appliquer l'inégalité de 
Sobolev du lemme~\ref{min:lem1}, ce qui garantit que la constante définie 
par le lemme est strictement positive.

Si on note $(g_u,\varphi_u)=(e^{2u}g,\varphi-\alpha u)$ un élément de 
la classe $[g,\varphi]_\alpha$, et si $\omega\in\Omega^p(M)$, alors
la norme $\|\omega\|_{r,\varphi_u}$ calculé pour la métrique
$g_u$ vaut:
\begin{eqnarray}
\int_M|\omega|_{g_u}^re^{-r\varphi_u}\de v_{g_u} & = &
\int_M e^{-p\frac n{p-\alpha}u}|\omega|_g^r
e^{-r\varphi-\frac n{p-\alpha}\alpha u}e^{nu}\de v_g\nonumber\\
& = & \int_M|\omega|_g^re^{-r\varphi}\de v_g.
\end{eqnarray}
Cette norme est donc constante sur la classe conforme pondérée. Par
conséquent, il en va de même pour le deux normes 
$\|\omega-\theta\|_{\frac n{p-\alpha},\varphi}$ et 
$\|\de\omega\|_{\frac n{p+1-\alpha},\varphi}$ qui définissent la 
constante de Sobolev et de cette constante elle-même.
\end{demo}

\begin{demo}[du théorème \ref{min:theo}]
Par commodité, on supposera dans la démonstration que la variété
$M^n$ est de volume 1.

Selon la proposition~\ref{theo:spectre3} on doit minorer le quotient
$\|\de\omega\|_{2,\varphi}^2/(\inf_{\de\theta=0}
\|\omega-\theta\|_{2,\varphi}^2)$.
Pour ce faire, on remarque la condition sur $\alpha$ peut se réécrire
$\frac n{p-\alpha}\geq2$ et $\frac n{p+1-\alpha}\leq2$. Donc, selon 
l'inégalité de Hölder, il existe des constantes $C_1,C_2>0$ dépendant de 
$g$ et $\varphi$ mais pas des formes $\omega$ et $\theta$ telles que 
$\|\omega-\theta\|_{2,\varphi}\leq 
C_1\|\omega-\theta\|_{\frac n{p-\alpha},\varphi}$ et 
$\|\de\omega\|_{\frac n{p+1-\alpha},\varphi}\leq 
C_2\|\de\omega\|_{2,\varphi}$. On en déduit
\begin{eqnarray}
\|\omega-\theta\|_{2,\varphi} & \leq &
C_1\|\omega-\theta\|_{\frac n{p-\alpha},\varphi}\leq C_1
K_p(M,[g,\varphi]_\alpha)\|\de\omega\|_{\frac n{p+1-\alpha},\varphi}\nonumber
\\
& \leq & C_1C_2K_p(M,[g,\varphi]_\alpha)\|\de\omega\|_{2,\varphi}
\end{eqnarray}
et donc
\begin{equation}
\frac{\|\de\omega\|_{2,\varphi}^2}{\inf_{\de\theta=0}
\|\omega-\theta\|_{2,\varphi}^2}\geq\frac1{C_1^2C_2^2
K_p(M,[g,\varphi]_\alpha)^2}.
\end{equation}
La constante $c$ du théorème est donc le membre de droite de l'inégalité
précédente.
\end{demo}

\section{Convergence du spectre vers celui d'un domaine}\label{conv}
Nous allons maintenant montrer qu'on peut faire tendre le spectre
du laplacien de Witten sur une variété compacte vers celui d'un
domaine. Ce résultat de convergence spectrale généralise à la fois
les théorèmes obtenus pour les laplaciens agissant sur
les fonctions (\cite{cdv86}) et sur le formes différentielles (\cite{ja11}),
et un résultat analogue obtenu par Y.~Colin de Verdière dans \cite{cdv87} 
pour le laplacien de Witten resteint aux fonctions.

Pour appliquer ce résultat, nous aurons besoin d'une
certaine uniformité de la convergence, nous reprendrons pour cela les
notations de \cite{cdv86} :

 Soit $E_0$ et $E_1$ sont deux sous-espaces vectoriels de même dimension~$N$
d'un espace de Hilbert, munis respectivement des formes quadratiques $q_0$
et $q_1$. Si $E_0$ et $E_1$ sont suffisamment proches, il existe
une isométrie naturelle $\psi$ entre les deux (voir la section~I de
\cite{cdv86} pour les détails de la construction), on définit alors
l'écart entre $q_0$ et $q_1$ par $\|q_1\circ\psi-q_0\|$. Pour deux
formes quadratiques $Q_0$ et $Q_1$ sur l'espace de Hilbert, on appellera
\emph{$N$-écart spectral entre $Q_0$ et $Q_1$} l'écart entre les
deux formes quadratiques restreintes à la somme des espaces propres associés
aux $N$ premières valeurs propres. Si cet écart est petit, alors les $N$
premières valeurs propres de $Q_0$ et leurs espaces propres sont proches
de ceux de $Q_1$.

 On veut montrer que la convergence spectrale est uniforme pour une certaine
famille de spectres limites. Comme dans \cite{cdv86} on dira donc qu'une forme
quadratique vérifie l'hypothèse ($*$) si ses valeurs propres vérifient
$$\lambda_1\leq\ldots\leq\lambda_N<\lambda_N+\eta\leq\lambda_{N+1}\leq M$$
pour un entier $N$ et des réels $\eta,M>0$ fixés une fois pour toute.

 Comme dans le cas du laplacien de Hodge-de~Rham, le procédé de convergence
fait apparaître un certain nombre de petites valeurs propres. Pour
les compter précisément, on introduit un espace de cohomologie qui traduit 
l'interaction entre la cohomologie du domaine $U$ et celle de la variété $M$,
et qui est défini comme le quotient des formes fermées de $U$ par la 
restriction des formes fermées de $M$ :
\begin{equation}
H^p(U/M)=\{\omega\in\Omega^p(\overline U),\ \de\omega=0\}/
\{\omega_{|\overline U},\ \omega\in\Omega^p(M)\textrm{ et }\de\omega=0\}
\end{equation}
Comme on l'a remarqué dans \cite{ja11}, cet espace est isomorphe au quotient 
de $H^p(U)$ par l'image de l'application naturelle $H^p(M)\to H^p(U)$ 
définie par restriction des formes fermées et exactes. En particulier, 
$H^p(U/M)$ est de dimension finie.

Le théorème que nous allons démontrer donne une convergence pour tous les 
degrés sauf un ou 
deux, en fixant la classe conforme pondérée :
\begin{theo}\label{conv:th2}
Soit $(M^n,g)$ une variété riemanienne compacte sans bord de dimension $n$
et $U$ un domaine de $M$ à bord $C^1$ dont le bord ne rencontre pas celui
de $M$, $\varphi$ une fonction sur $M$ et $\alpha>0$ un réel. 
Il existe une suite de métriques $g_i$ et une suite de fonctions 
$\varphi_i$ sur $M$ et une constante $c>0$ telles que pour tout $i$, 
$(g_i,\varphi_i)\in[g,\varphi]_\alpha$ et
\begin{enumerate}
\item $\tilde\mu_{p,k}(M,g_i,\varphi_i)\to0$ pour tout
$k\leq d_p$ et tout $p<\frac n2+\alpha-1$ quand $i\to\infty$;
\item $\tilde\mu_{p,k+d_p}(M,g_i,\varphi_i)\to\tilde\mu_{p,k}(U,g,\varphi)$
pour tout $k\geq1$ et tout $p<\frac n2+\alpha-1$ quand $i\to\infty$;
\item $\tilde\mu_{p,k}(M,g_i,\varphi_i)\Vol(M,g)^{\frac2n}>c$ pour 
$p\in[\frac n2+\alpha-1,\frac n2+\alpha]$;
\item $\tilde\mu_{p,k}(M,g_i,\varphi_i)\to0$ pour tout
$k\leq d_{n-p-1}$ et tout $p\in]\frac n2+\alpha,n-1]$ quand $i\to\infty$;
\item $\tilde\mu_{p,k+d_p}(M,g_i,\varphi_i)\to\tilde\mu_{n-p-1,k}(U,g,-\varphi)$
pour tout $k\geq1$ et tout $p\in]\frac n2+\alpha,n-1]$ quand $i\to\infty$;
\end{enumerate}
où $d_p$ est la dimension de $H^p(U/M)$.

En outre, si les $\tilde\mu_{p,k}(U,g,\varphi)$ vérifient l'hypothèse ($*$) 
pour un $p$ donné, alors pour tout $\varepsilon>0$ il existe
$i$ tel que le $N$-écart spectral entre les laplaciens de Witten sur
$\tilde\Delta_\varphi$ sur $U$ et $\tilde\Delta_{\varphi_{i}}$ sur $M$
soit inférieur à $\varepsilon$.
\end{theo}

 Grâce à la proposition~\ref{theo:spectre3}, la démonstration est
très similaire à celle du théorème~11 de \cite{ja11}. Cependant, les
modifications et précisions à apporter, bien que mineures, sont assez
nombreuses. Nous allons donc donner ici la démonstration complète.

Commençons par rappeler deux lemmes de~\cite{cdv86} sur lesquels nous
nous appuieront.  Les constantes $N$, $M$ et $\eta$ qui interviennent 
dans les énoncés font référence à l'hypothèse ($*$) définie plus haut.

\begin{lem}[\cite{cdv86}, th.~I.7]\label{conv:lem1}
Soit $Q$ une forme quadratique positive sur un espace de Hilbert $\mathcal H$
dont le domaine admet la décomposition $Q$-orthogonale
$\mathrm{dom}(Q)=\mathcal H_0\oplus\mathcal H_\infty$. Pour tout
$\varepsilon>0$, il existe une constante $C(\eta,M,N,\varepsilon)>0$ (grande)
telle que si $Q_0=Q_{|\mathcal H_0}$ vérifie l'hypothèse ($*$) et que
$\forall x\in\mathcal H_\infty,\ Q(x)\geq C|x|^2$,
alors $Q$ et $Q_0$ ont un $N$-écart spectral inférieur à $\varepsilon$.
\end{lem}

\begin{lem}[\cite{cdv86}, th.~I.8]\label{conv:lem2}
Soit $(\mathcal H,|\cdot|)$ un espace de Hilbert muni d'une forme quadratique
positive $Q$. On se donne en outre une suite de métriques $|\cdot|_n$
sur $\mathcal H$ et une suite de formes quadratiques $Q_n$ de même
domaine que $Q$ telles que:
\begin{itemize}
\item[(i)] il existe $C_1,C_2>0$ tels que $\forall x\in\mathcal H,\
C_1|x|\leq |x|_n\leq C_2|x|$;
\item[(ii)] pour tout $x\in \mathrm{dom}(Q)$, $|x|_n\to|x|$;
\item[(iii)] pour tout $x\in \mathrm{dom}(Q)$, $Q(x)\leq Q_n(x)$;
\item[(iv)] pour tout $x\in \mathrm{dom}(Q)$, $Q_n(x)\to Q(x)$.
\end{itemize}
Si $Q$ vérifie l'hypothèse ($*$), alors à partir d'un certain rang (dépendant
de $\eta$, $M$ et $N$), $Q$ et $Q_n$ ont un $N$-écart spectral inférieur à
$\varepsilon$.
\end{lem}
\begin{remarque}\label{conv:rem1}
Comme on l'a remarqué dans \cite{ja11}, dans le lemme~\ref{conv:lem2}, on peut 
affaiblir l'hypothèse~(i) en $C_1|x|\leq |x|_n\leq C_2|x|+
\varepsilon_nQ(x)^\frac12$ avec $\varepsilon_n\to0$,
la démonstration restant exactement la même.
\end{remarque}

D'après la proposition~\ref{theo:spectre3}, on veut montrer la convergence
de la forme quadratique $Q(\omega)=\int_M |\omega|^2\de v_\varphi$ par rapport
à la norme $|\omega|=\inf_{\de\theta=\omega}
\int_M|\theta|^2\de v_\varphi$, pour tout degré $p$.

On passe par l'intermédiaire d'une famille de métriques et de 
fonctions singulières $(g_\varepsilon,\varphi_\varepsilon)$, 
$\varepsilon\in]0,1]$ définie par :
\begin{equation}
\left\{\begin{array}{ccl}
(g_\varepsilon,\varphi_\varepsilon) & = & (g,\varphi) \textrm{ sur } U,\\
(g_\varepsilon,\varphi_\varepsilon) & = & (\varepsilon^2g,\varphi
-\alpha\ln\varepsilon)\textrm{ sur } M\backslash U.
\end{array}\right.
\end{equation}
Le couple $(g_\varepsilon,\varphi_\varepsilon)$ est bien contenu dans
la classe conforme (singulière) pondérée $[g,\varphi]_\alpha$.
 La forme quadratique $Q_\varepsilon$ et
la norme $|\cdot|_\varepsilon$ associées sont alors données,
pour une forme $\omega\in\Omega^{p+1}(M)$, par
\begin{equation}\label{conv:eq1}
Q_\varepsilon(\omega)=\int_U|\omega|^2e^{-2\varphi}\de v_g
+\varepsilon^{n+2\alpha-2p-2}
\int_{M\backslash U}|\omega|^2e^{-2\varphi}\de v_g
\end{equation}
et
\begin{equation}\label{conv:eq2}
|\omega|_\varepsilon^2=\inf_{\de\theta=\omega}\left(\int_U|\theta|^2
e^{-2\varphi}\de v_g+\varepsilon^{n+2\alpha-2p}
\int_{M\backslash U}|\theta|^2e^{-2\varphi}\de v_g\right).
\end{equation}
On voit que, sous la condition $p<\frac n2+\alpha-1$, le dernier
terme de chacune des deux expressions précédentes tend vers 0.
\begin{lem}\label{conv:lissage}
Pour tout $\varepsilon>0$, il existe une suite $(g_j,\varphi_j)\in
[g,\varphi]_\alpha$ formée
de métriques et de fonctions lisses tendant vers 
$(g_\varepsilon,\varphi_\varepsilon)$ et telle que le volume de $(M,g_j)$
tende vers celui de $(M,g_\varepsilon)$ et pour tout $p<\frac n2+\alpha-1$
et que pour tout $k\geq1$, la valeur propre $\tilde\mu_{p,k}(M,g_j,\varphi_j)$
tende vers $\tilde\mu_{p,k}(M,g_\varepsilon,\varphi_\varepsilon)$, avec 
convergence des espaces propres.
\end{lem}
\begin{demo}
Le principe de la démonstration consiste à chercher à appliquer le
lemme~\ref{conv:lem2}.

Pour un $\varepsilon>0$ fixé, on peut approcher la fonction
$\chi_U+\varepsilon\chi_{M\backslash U}$ par une suite de fonctions
décroissantes $(f_j)$ vérifiant $f_j\leq1$. On définit alors la suite 
$(g_j,\varphi_j)$ par
$g_j=f_j^2\cdot g$ et $\varphi_j=\varphi-\alpha\ln f_j$ et 
on note $Q_j$ et $|\cdot|_j$ la forme quadratique et la norme 
hilbertienne associées. On a bien $(g_j,\varphi_j)\in
[g,\varphi]_\alpha$, la convergence du volume est immédiate, et on 
peut écrire
$Q_j(\omega)=\int_Mf_j^{n+2\alpha-2p-2}|\omega|^2e^{-2\varphi}\de v_g$ et
$|\omega|_j^2=\inf_{\de\theta=\omega}(\int_Mf_j^{n+2\alpha-2p}
|\theta|^2e^{-2\varphi}\de v_g)$.

Comme $p<\frac n2+\alpha-1$, on peut vérifier que la famille de forme 
quadratique $Q_j$ est décroissante et converge simplement vers 
$Q_\varepsilon$, c'est-à-dire que les hypothèses (iii) et (iv) du 
lemme~\ref{conv:lem2} sont satisfaites.

Comme la suite $(f_j)$ vérifie $\chi_U+\varepsilon\chi_{M\backslash U}
\leq(f_j)\leq\frac1\varepsilon(\chi_U+\varepsilon\chi_{M\backslash U})$,
on a $|\omega|_\varepsilon^2\leq |\omega|_j^2\leq \frac1\varepsilon\cdot 
|\omega|_\varepsilon^2$ pour toute forme exacte $\omega$, ce qui constitue
l'hypothèse~(i) du lemme~\ref{conv:lem2}.
 
Il reste à vérifier l'hypothèse (ii), c'est-a-dire que $|\cdot|_j$
converge simplement vers $|\cdot|_\varepsilon$. Soit $\eta>0$. Par
définition de $|\cdot|_\varepsilon$, il existe une forme $\theta_0$ telle 
que $\de\theta_0=\omega$ et $\|\theta_0\|_{g_\varepsilon,\varphi}^2\leq
|\omega|_\varepsilon^2+\eta$.
Comme $(g_j,\varphi_j)$ converge vers $(g_\varepsilon,\varphi_\varepsilon)$, 
pour $j$ assez grand on a aussi $\|\theta_0\|_{g_j,\varphi_j}^2\leq
\|\theta_0\|_{g_\varepsilon,\varphi_\varepsilon}^2+\eta$  et donc aussi
$|\omega|_\varepsilon^2\leq|\omega|_j^2\leq 
\|\theta_0\|_{g_j,\varphi_\varepsilon}^2\leq |\omega|_\varepsilon^2+2\eta$. 
Par conséquent on a bien $|\omega|_j\to|\omega|_\varepsilon$ pour tout 
$\omega$.

 Selon le lemme~\ref{conv:lem2}, à $\varepsilon$ fixé, on peut donc trouver
un $j_\varepsilon$ tel que l'écart spectral entre $Q_{j_\varepsilon}$ et
$Q_\varepsilon$ soit arbitrairement petit.
\end{demo}

\begin{demo}[du théorème~\ref{conv:th2}]
Remarquons d'abord que le point~3 du théorème, c'est-à-dire le cas
$\frac n2+\alpha-1\leq p\leq\frac n2+\alpha$, découle de la minoration 
conforme donnée par le théorème~\ref{min:theo}.

Par ailleurs, le cas $p>\frac n2+\alpha$ se déduit du cas 
$p<\frac n2+\alpha-1$. En effet, la relation de commutation
(\ref{theo:hodge}) avec la dualité de Hodge donne 
$\tilde\mu_{p,k}(g_i,\varphi_i)=\tilde\mu_{n-p-1,k}(g_i,-\varphi_i)$.

Il reste donc a traiter le cas $p<\frac n2+\alpha-1$. Grâce au 
lemme~\ref{conv:lissage}, on est ramené à montrer la convergence du
spectre de $Q_\varepsilon$. La démonstration se déroule en deux 
étapes. La première consiste à décomposer l'espace
$\mathcal H$ des formes exactes en une somme $\mathcal H_0\oplus
\mathcal H_\infty$ à laquelle on applique le lemme~\ref{conv:lem1}.
Dans la seconde, on montre la convergence du spectre de $Q_\varepsilon$ 
restreint à $\mathcal H_0$ vers le spectre du domaine à l'aide du 
lemme~\ref{conv:lem2}. Sauf mention contraire, les normes $L^2$ seront 
considérées relativement à la mesure $\de v_\varphi=e^{-2\varphi}\de v_g$.

\emph{Étape~1.}
On commence par définir le sous-espace $\mathcal H_\infty$ de
$\overline{\Ima\de^p}\subset L^2(\Lambda^{p+1}M)$ comme l'adhérence des
différentielles des formes lisses qui s'annulent sur $U$ (les adhérences
sont au sens de la norme $L^2$). Une telle forme
va nécessairement vérifier la condition de bord de Dirichlet sur
$M\backslash U$:
\begin{equation}
\mathcal H_\infty=\overline{\{\de\theta,\ \theta\in\Omega^p(M),\ 
\theta_{|U}=0,\ \theta_{|M\backslash U}\textrm{ vérifie }(D)\}}.
\end{equation}

L'espace $\mathcal H_0$ est défini comme la somme de deux espaces
$\mathcal H_1$ et $\mathcal H_2$ construits séparément.

 Soit $\omega$ une $(p+1)$-forme exacte sur $\overline U$ et
$\theta\in\Omega^p(\overline U)$ la primitive de $\omega$ qui minimise
la norme $L^2(U,\de v_\varphi)$. Cette forme
$\theta$ est obtenue par application du lemme~\ref{theo:prim} et
de la remarque~\ref{theo:rem}. Si $\tilde\theta$ est un prolongement lisse
de $\varphi$ sur $M$, alors $\de\tilde\theta$ est définie à un élément
de $\mathcal H_\infty$ près. On peut définir alors $\tilde\omega$ comme
le $\de\tilde\theta$ de norme minimale pour $g_\varepsilon$. Cet infimum
est bien atteint dans $L^2(M,\de v_\varphi)$ et on peut le construire par
projection sur l'orthogonal $L^2$ de $\mathcal H_\infty$. On pose alors
\begin{equation}
\mathcal H_1=\overline{\{\tilde\omega,\ \omega\in\Omega^{p+1}(\overline U)
\textrm{ exacte}\}}.
\end{equation}

 L'espace $\mathcal H_2$ est défini à partir de l'espace de cohomologie
$H^p(U/M)$. Comme $H^p(U/M)$ est isomorphe au quotient de $H^p(U)$ par le
sous-espace induit par $H^p(M)$, il est aussi isomorphe à un supplémentaire,
dans l'espace des formes harmoniques (pour le laplacien de Witten) de 
$U$ avec condition de bord ($A_\varphi$),
des représentants harmoniques des classes de cohomologie induites par
$H^p(M)$. Grâce au  corollaire~\ref{theo:coh}, chaque classe 
$[c]\in H^p(U/M)$ possède un représentant $T_\varphi h$ où $h$ est 
$\tilde\Delta_\varphi$-harmonique et tel que $T_\varphi h$ soit
$L^2(U,\de v_\varphi)$-orthogonale aux restrictions des formes fermées de $M$.

On peut alors construire $\mathcal H_2$ sur le modèle de $\mathcal H_1$.
Chaque forme $T_\varphi h$ représentant une classe de $H^p(U/M)$
peut être étendue en une forme $\tilde h$ sur $M$, la forme $\de\tilde h$
étant alors définie à un élément de $\mathcal H_\infty$ près. On notant
$\tilde\omega_h$ la forme $\de\tilde h$ qui est 
$L^2(M,\de v_\varphi)$-orthogonale à $\mathcal H_\infty$, et on pose
\begin{equation}
\mathcal H_2=\{\tilde\omega_h,\ [h]\in H^p(U/M)\}.
\end{equation}
 
Si on pose $\mathcal H_0=\mathcal H_1\oplus\mathcal H_2$, alors
$\mathcal H_0$ et  $\mathcal H_\infty$ sont par construction orthogonaux
pour la norme $L^2(U,\de v_\varphi)$, donc $Q_\varepsilon$-orthogonaux.
Avant d'appliquer le lemme~\ref{conv:lem1} on doit encore vérifier que
$\mathcal H_0\oplus\mathcal H_\infty$ contient bien toutes les
$(p+1)$-formes exactes. Soit $\omega\in\Omega^{p+1}$ une forme exacte. Par
définition de $\mathcal H_1$, on peut écrire $\omega=\omega_1+\omega'$, avec
$\omega_1\in\mathcal H_1$ telle que $\omega_{|U}={\omega_1}_{|U}$ et
$\omega'_{|U}=0$. Comme $\omega$ et $\omega_1$
sont exacte, $\omega'$ l'est aussi. Si on pose $\omega'=\de\theta$, la
forme $\theta$ est définie à une forme fermée près et vérifie
$\de\theta=0$ sur $U$, la classe $[\theta]\in H^p(U/M)$ est donc
bien définie. Par définition de $\mathcal H_2$ on a alors $\omega'=
\tilde\omega_h+\omega_0$, où $\tilde\omega_h$ est l'élément de $\mathcal H_2$
associé au représentant harmonique $h$ de $[\theta]\in H^p(U/M)$, et
$\omega_0\in\mathcal H_\infty$. On a donc bien $\omega\in\mathcal H_1
\oplus\mathcal H_2\oplus\mathcal H_\infty$, et donc
$\mathrm{dom}(Q_\varepsilon)=\mathcal H_0\oplus\mathcal H_\infty$.

Si on note $\lambda^{(D)}$ la première valeur propre du laplacien 
de Hodge sur $M\backslash U$ pour la métrique $g$ et la condition de bord (D), 
on sait que pour toute forme $\omega\in \mathcal H_\infty$ il existe, par 
définition, une forme $\theta\in\Omega^p(M)$ à support dans $M\backslash U$ 
telle que $\de\theta=\omega$ et $\int_{M\backslash U}|\omega|^2/
\int_{M\backslash U}|\theta|^2\geq\lambda^{(D)}$ 
pour la métrique $g$. Il existe donc une constante $C_\varphi$ ne dépendant
que de $\varphi$ telle que $\|\omega\|^2/\|\theta\|^2\geq C_\varphi
\lambda^{(D)}$.  Pour la métrique $g_\varepsilon$, on a alors
$\|\omega\|^2/\|\theta\|^2\geq\varepsilon^{-2}C_\varphi\lambda^{(D)}$, et 
\emph{a fortiori} $Q_\varepsilon(\omega)/|\omega|_\varepsilon\geq
\varepsilon^{-2}C_\varphi\lambda^{(D)}$. Si $\varepsilon$ est suffisamment 
petit, on peut appliquer le lemme~\ref{conv:lem1}.

\emph{Étape~2.}
On va maintenant chercher à appliquer le lemme~\ref{conv:lem2} et
la remarque~\ref{conv:rem1} à l'espace $\mathcal H_0$ et aux familles
de métriques et de formes quadratiques $|\cdot|_\varepsilon$ et
$Q_\varepsilon$.  On définit la forme quadratique $Q$ sur $\mathcal H_0$
par $Q(\omega)=\int_U|\omega|^2\de v_\varphi$. Pour définir une
norme $|\cdot|$ sur $\mathcal H_0$, on procède de la manière suivante : 
pour $\omega\in\mathcal H_1$, on note $\theta_\omega$ l'image par $T_\varphi$
de la primitive coexacte de $\omega_{|\overline U}$ donnée par la
proposition~\ref{theo:prim}, et pour $\omega\in\mathcal H_2$, on note
$\theta_\omega$ l'image par $T_\varphi$ du représentant harmonique de la 
classe de cohomologie définie par $\omega$. On étend linéairement 
l'application $\omega\to\theta_\omega$ et on pose 
$|\omega|=\|\theta_\omega\|$, la norme $\|\cdot\|$ étant ici la norme 
$L^2(U,\de v_\varphi)$ sur les $p$-formes de~$U$.
Les espaces $\mathcal H_1$ et $\mathcal H_2$ sont orthogonaux pour
$\|\cdot\|$, le noyau de la forme quadratique $Q$ (relativement à $\|\cdot\|$)
est $\mathcal H_2$, et le spectre de $Q$ sur $\mathcal H_1$ est le spectre
du domaine $U$ pour les $(p+1)$-formes exactes.

 On doit maintenant vérifier que les quatre hypothèses du
lemme~\ref{conv:lem2} sont satisfaites.

Les hypothèses~(iii) et~(iv) sont les plus simples à vérifier~: par
définition de $Q_\varepsilon$, pour tout $\omega\in \mathcal H_0$, 
$Q_\varepsilon(\omega)$ tend vers $Q(\omega)=\int_U|\omega|^2\de v_\varphi$ 
et $Q_\varepsilon\geq Q$ pour tout $\varepsilon$.

Passons à l'hypothèse (ii).
 Par définition de $\mathcal H_1$ et $\mathcal H_2$, il existe
pour tout $\omega\in\mathcal H_0$ un prolongement $\tilde\theta_\omega$ 
de $\theta_\omega$ tel que $\de\tilde\theta_\omega=\omega$.
On a alors
\begin{equation}
|\omega|_\varepsilon^2\leq\|\tilde\theta_\omega\|_\varepsilon^2\to
\int_U|\tilde\theta_\omega|^2\de v_\varphi=|\omega|^2.
\end{equation}
En outre, si $\de\theta=\omega$, alors $\theta$ et $\tilde\theta_\omega$
ne diffèrent que par une forme fermée de $M$, et donc leur restriction à
$\overline U$ ne diffèrent aussi que par une forme fermée. Par conséquent, la
norme de $\theta_\omega$ minore la norme $L^2(U,\de v_\varphi)$ de $\theta$ 
(cf. corollaire~\ref{theo:coh} et théorème~\ref{theo:prim}), et
\begin{equation}
|\omega|^2=\int_U|\tilde\theta_\omega|^2\de v_g\leq\int_U|\theta|^2\de v_g
\leq\|\theta\|_\varepsilon^2
\end{equation}
On en déduit que $|\omega|\leq|\omega|_\varepsilon$, et donc
que $|\omega|_\varepsilon\to|\omega|$ pour tout $\omega$.

L'hypothèse (i), dans sa version faible (cf. remarque~\ref{conv:rem1}),
est la plus technique à vérifier.
On doit contrôler $|\cdot|_\varepsilon$ en fonction de $|\cdot|$ et
$Q$. Pour ce faire, on fixe un élément $\omega\in\mathcal H_0$ et on va
construire une primitive particulière dont la norme $L^2$ majorera
$|\omega|_\varepsilon$.

 On note $\bar\theta$ la $p$-forme définie par $\bar\theta=\theta_\omega$
sur $\overline U$, et prolongée harmoniquement (pour le laplacien de Hodge et
la métrique $g$) sur $M$.On a alors $\|\bar\theta\|_{H^1(M\backslash U)}\leq C 
\|\theta_{|\partial U}\|_{H^\frac12(\partial U)}$, la constante $C$ dépendant 
de $\varphi$ mais pas de $\theta$ (la démonstration dans le cas $\varphi=0$
est similaire au cas des fonctions, cf.~\cite{ta96}, ch.~5, prop.~1.7. Le
cas général découle du fait que $\varphi$ est bornée donc que les normes
pour $\de v_g$ et $\de v_\varphi$ sont équivalentes).

 Par définition de
$\mathcal H_0$, la $(p+1)$-forme $\bar\omega_\infty=\de\bar\theta-\omega$
est un élément de $\mathcal H_\infty$ et elle est orthogonale à $\omega$.
Elle vérifie donc $\|\bar\omega_\infty\|^2_{L^2(M\backslash U)}
\leq\|\de\bar\theta\|^2_{L^2(M\backslash U)}\leq
\|\bar\theta\|_{H^1(M\backslash U)}^2$ et
elle admet une primitive $\bar\theta_\infty$ nulle sur $U$ dont la norme
vérifie $\|\bar\theta_\infty\|^2\leq\|\bar\omega_\infty\|^2/
(C_\varphi\lambda^{(D)})$. Si on pose $\theta=\bar\theta-\bar\theta_\infty$,
on a alors $\de\theta=\omega$ et
\begin{eqnarray}
\|\theta\|_{L^2(M\backslash U)}&\leq&\|\bar\theta\|_{L^2(M\backslash U)}+
\|\bar\theta_\infty\|_{L^2(M\backslash U)}\nonumber\\
&\leq& \|\bar\theta\|_{L^2(M\backslash U)}+
\|\bar\theta\|_{H^1(M\backslash U)} /{(C_\varphi\lambda^{(D)})}^{1/2},
\end{eqnarray}
Toutes les normes étant ici relative à la métrique $g$ et la mesure 
$\de v_\varphi$.

Comme $\|\bar\theta\|_{H^1(M\backslash U)}\leq C 
\|\theta_{|\partial U}\|_{H^\frac12(\partial U)}$ et que la norme
$\|\theta_{|\partial U}\|_{H^\frac12(\partial U)}$ est elle-même contrôlée
par la norme $H^1$ de $\theta$ sur $U$, on a finalement
\begin{equation}\label{conv:eq3}
\int_{M\backslash U}|\theta|^2\de v_g\leq C'\|\theta\|_{H^1(U)},
\end{equation}
où $C'$ est une constante dépendant de $g$ et $\varphi$ mais pas de 
$\varepsilon$.

En utilisant le fait que $\theta$ est cofermée sur $U$ (car égale à
$\theta_\omega$) et tangentielle le long de $\partial U$, une inégalité
elliptique à trace associée à l'opérateur $\de+\codiff$ (voir~\cite{ta96},
section~5.9) donne, en utilisant encore une fois que les normes
associées à $\de v_g$ et $\de v_\varphi$ sont équivalentes,
\begin{equation}\label{conv:eq4}
\|\theta\|_{H^1(U)}
\leq C''(\|\theta\|_{L^2(U)}+\|\de\theta\|_{L^2(U)})=
C''(|\omega|+Q(\omega)^\frac12),
\end{equation}
la métrique considérée étant ici encore $g$.

Pour une métrique $g_\varepsilon$, on a $|\omega|_\varepsilon^2\leq
\|\theta\|_{L^2(U)}^2+\varepsilon^{n-2p}\|\theta\|_{L^2(M\backslash U)}^2$.
Comme $\|\theta\|_{L^2(U)}=|\omega|$ et que $\|\theta\|_{L^2(M\backslash U)}$
peut être majoré à l'aide de~(\ref{conv:eq3}) et~(\ref{conv:eq4}),
on obtient la majoration
\begin{equation}
|\omega|_\varepsilon^2\leq |\omega|^2+\varepsilon^{n-2p}C'C''^2
(|\omega|+Q(\omega)^\frac12)^2
\end{equation}
qui permet d'appliquer la remarque~\ref{conv:rem1} et le lemme~\ref{conv:lem2}.
\end{demo}

\section{Spectre d'une variété privée d'une boule}\label{boule}
Cette section est consacrée au théorème qui suit sur la convergence
du spectre sur une variété privée d'une boule :
\begin{theo}\label{boule:theo}
Soit $(M,g)$  une variété riemannienne compacte, $\varphi$ une fonction
lisse sur $M$ et $x\in M$. Pour tout $\varepsilon>0$, on note
$M_\varepsilon=M\backslash B(x,\varepsilon)$ la variété $M$ privée
d'une boule de rayon $\varepsilon$. Alors il existe une famille de 
fonction $\varphi_\varepsilon$ tendant vers $\varphi$ pour la 
topologie $C^0$ telle que
$\tilde\mu_{p,k}(M_\varepsilon,g,\varphi_\varepsilon)\to
\tilde\mu_{p,k}(M,g,\varphi)$
pour tout $p\leq n-1$ et $i>0$, avec convergence des espaces propres.
\end{theo}
Dans le cas du laplacien de Hodge-de~Rham, c'est-à-dire avec 
$\varphi_\varepsilon=\varphi=0$, ce théorème a été montré
par C.~Anné et B.~Colbois dans \cite{ac93} (théorème~1.1). La comparaisons 
des deux énoncés imposent deux précisions. D'une part, dans \cite{ac93},
la convergence du spectre des $p$-formes est montrée sans distinguer 
formes exactes et coexactes, mais par théorie de Hodge cela implique 
la convergence en restriction aux formes coexactes. D'autre part, le
résultat de \cite{ac93} fait apparaître une valeur propre tendant
vers~0 pour les $n-1$ formes, mais la forme propre est nécessairement
exacte et la petite valeur propre correspond à celle qui apparaît 
pour les fonctions avec condition de Dirichlet.
 
La démonstration de \cite{ac93} fonctionne correctement pour le laplacien 
de Witten si, sur une boule de rayon $R$ fixé et centrée en $x$, la fonction 
$\varphi$ est constante. On a alors $\tilde\Delta_\varphi=\Delta$, et en 
particulier le lemme suivant 
(\cite{ac93}, corollaire de la proposition~3.3), qui concentre l'essentiel 
des difficultés, s'applique parfaitement sous cette hypothèse:
\begin{lem}[\cite{ac93}]
Il existe une constante $C>0$ telle que pour toute forme $\omega
\in H^1(\Lambda^p(B(x,R)-B(x,\varepsilon)))$ tangentielle le long de
$B(x,\varepsilon)$ et nulle le long de $B(x,R)$, on a
$$\|P\omega\|_q\leq C\|\omega\|_q,$$
où $P\omega$ désigne le prolongement harmonique de $\omega$ sur 
$B(x,R)$ et $\|\cdot\|$ la norme d'opérateur associé au laplacien de 
Hodge-de~Rham.
\end{lem}

On va donc montrer qu'on peut se ramener à cette situation:
\begin{lem}\label{boule:lem}
Soit $(M,g)$  une variété riemannienne compacte, $\varphi$ une fonction
lisse sur $M$ et $x\in M$. Il existe une suite de fonction $\varphi_i$ telle 
que
\begin{enumerate}
\item $\varphi_i$ est constante sur $B(x,\frac1i)$;
\item $\varphi_i\to\varphi$ pour la topologie $C^0$;
\item $\tilde\mu_{p,k}(M,g,\varphi_i)\to\tilde\mu_{p,k}(M,g,\varphi_i)$
avec convergence des espaces propres.
\end{enumerate}
\end{lem}
Pour approcher le spectre de $(M,g,\varphi)$, on commence donc par
l'approcher par celui de $(M,g,\varphi_i)$ à l'aide de ce lemme, 
puis on applique les résultats de \cite{ac93}.

\begin{demo}[du lemme~\ref{boule:lem}]
On peut déformer la fonction $\varphi$ en une suite $\varphi_i$ telle que 
$\varphi_i$ soit constante sur $B(x,\frac1i)$, et de manière à ce que 
$\varphi_i$ converge vers $\varphi$ pour la topologie $C^0$.
En notant $Q$ et $|\cdot|$ (resp. $Q_i$ et $|\cdot|_i$) la forme quadratique
et la norme associés à $(g,\varphi)$ (resp. $(g,\varphi_i)$) comme dans
la proposition~\ref{theo:spectre3},
on a alors une suite décroissante $C_i$ telle que $C_i\to1$ et
\begin{equation}
\left\{\begin{array}{l}
\displaystyle\frac1{C_i}|\omega|\leq|\omega|_i\leq C_i|\omega|\\
\displaystyle\frac1{C_i}Q(\omega)\leq Q_i(\omega)\leq C_iQ(\omega)
\end{array}\right.
\end{equation}
Pour tout forme exacte $\omega$. En outre, on peut ajuster $\varphi_i$ de 
manière à avoir de plus $Q(\omega)\leq Q_i(\omega)$. On peut alors appliquer 
le lemme~\ref{conv:lem2}.
\end{demo}

\section{Application à la multiplicité}\label{multi}
Pour prescrire la multiplicité des valeurs propres du laplacien de Witten
nous utilisons, selon la méthode introduite par Colin de Verdière, trois
ingrédients : les théorèmes de convergence spectrale démontrés dans
les sections précédentes, des modèles de valeurs propres multiples déjà
connus et une propriété de stabilité vérifiée par ces modèles. Nous 
allons commencer par rappeler cette dernière.

On suppose qu'on a une famille d'opérateurs $(P_a)_{a\in B^k}$, où $B^k$
est la boule unité de $\R^k$ (en pratique, $P_a$ est le laplacien de Witten
associé à un couple $(g_a,\varphi_a)$), tels que
$P_0$ possède une valeur propre $\lambda_0$ d'espace propre $E_0$ et
de multiplicité $N$. Pour les petites valeurs de $a$, $P_a$ possède des
valeurs propres proches de $\lambda_0$ dont la somme des espaces propres
est de dimension~$N$. Comme dans la définition de l'écart spectral,
on identifie cette somme à $E_0$ et on note $q_a$ la forme
quadratique associée à $P_a$ transportée sur $E_0$.
\begin{definition}[\cite{cdv88}]\label{presc:def}
On dit que $\lambda_0$ vérifie l'hypothèse de transversalité d'Arnol'd si
l'application $\Psi:a\mapsto q_a$ de $B^k$ dans $\mathcal Q(E_0)$ est
essentielle en $0$, c'est-à-dire qu'il existe $\varepsilon>0$ tel que si
$\Phi:B^k\to\mathcal Q(E_0)$ vérifie $\|\Psi-\Phi\|_{\infty}\leq\varepsilon$,
alors il existe $a_0\in B^k$ tel que $\Phi(a_0)=q_0$.
\end{definition}
Une propriété cruciale est que si $\Phi$ provient d'une famille $(P'_a)$
d'opérateurs, alors $\lambda_0$ est valeur propre de $P'_{a_0}$ de
multiplicité~$N$ et vérifie la même propriété de transversalité, ce
qui justifie qu'on parle de stabilité de la multiplicité.  Comme remarqué 
dans \cite{cdv88}, on peut généraliser cette définition à une suite finie 
de valeurs propres.

Cette propriété est aussi préservée par produit, grâce à la forme 
de Künneth (théorème~\ref{theo:kunneth}):
\begin{lem}\label{multi:produit}
On suppose que $\tilde\mu_{p_1,k_1}(M_1,g_1,\varphi_1)$ est une valeur propre
simple, $\tilde\mu_{p_2,k_2}(M_2,g_2,\varphi_2)$ une valeur propre de 
multiplicité stable $N$ et qu'il n'y a pas d'autres valeurs propres de $M_1$ 
et $M_2$ dont la somme soit $\tilde\mu_{p_1,k_1}(M_1,g_1,\varphi_1)+
\tilde\mu_{p_2,k_2}(M_2,g_2,\varphi_2)$.
Alors $\tilde\mu_{p_1,k_1}(M_1,g_1,\varphi_1)+
\tilde\mu_{p_2,k_2}(M_2,g_2,\varphi_2)$ est valeur propre de multiplicité 
stable $N$ sur $(M_1\times M_2,g_1\oplus g_2,\varphi_1+\varphi_2)$.
\end{lem}
La démonstration est la même que dans le cas du laplacien de Hodge 
(\cite{ja11}, lemme~3.3). On utilisera en particulier ce lemme dans le cas
où $(M_2,g_2)$ est un intervalle  $[-\varepsilon,\varepsilon]$
avec $\varphi_2=0$. En effet, le noyau du laplacien de Hodge sur l'intervalle
pour la condition de bord (A) est de dimension~1 en degré~0 et trivial en
degré~1 (les fonctions harmoniques étant 
constantes), et toutes les autres valeurs propres tendent vers l'infini 
quand $\varepsilon\to0$. Sur $M_1\times [-\varepsilon,\varepsilon]$ les 
premières formes propres sont les relevés des formes propres de 
$(M_1,g_1,\varphi_1)$.

 Comme dans \cite{ja11}, la démarche pour démontrer les 
théorèmes~\ref{intro:th1} et \ref{intro:th2} consiste à d'abord construire,
pour tout $p$ et $n$, un modèle de valeur propre multiple pour les $p$-formes 
sur une variété de dimension~$n$, puis à obtenir plusieurs valeurs propres 
multiples sur une même variété en plongeant plusieurs modèles et en
appliquant le théorème de convergence~\ref{conv:th2}.

Commençons par rappeler un modèle de multiplicité qui servira de
point de départ pour la suite. On note $\mu_{p,i}$ les valeurs propres
du laplacien de Hodge, c'est-à-dire que $\mu_{p,i}(M,g)=
\tilde\mu_{p,i}(M,g,0)$. 
\begin{lem}[\cite{ja11}, lemme~3.4]\label{multi:lem1}
Pour tous entiers $N\geq1$, $n\geq3$, toute suite finie $0<a_1\leq a_2\leq
\ldots\leq a_N$ et toute constante $C>a_N$, il existe une métrique $g$ sur
$S^n$ telle que $\mu_{0,i}(S^n,g)=a_i$ pour $i\leq N$, ces valeurs propres
vérifiant l'hypothèse de stabilité, $\mu_{0,N+1}(S^n,g)>C$ et
$\mu_{p,1}(S^n,g)>C$ pour $1\leq p\leq[\frac{n-1}2]$, le volume $\Vol(S^n,g)$ 
étant arbitrairement petit.  
\end{lem}

Comme dans \cite{ja11}, on en déduit des modèles de valeurs propres 
multiples pour les autres degrés :
\begin{lem}\label{multi:lem2}
Pour tous entiers $N\geq1$, $p\geq2$ et $n\geq p+1$, toute suite finie
$0<a_1\leq a_2\leq\ldots\leq a_N$ et toute constante $C>a_N$, il existe
une métrique $g$ et une fonction $\varphi$ sur la variété 
$S^n$ telles que $\tilde\mu_{p,i}(S^n,g,\varphi)=a_i$ pour
$i\leq N$, ces valeurs propres vérifiant l'hypothèse de stabilité,
$\tilde\mu_{p,N+1}(S^n,g,\varphi)>C$ et $\tilde\mu_{q,1}(S^n,g,\varphi)>C$ 
pour $1\leq q\leq n-1$, $q\neq p$, le volume $\Vol(S^n,g)$ étant 
arbitrairement petit.
\end{lem}
\begin{demo}
On procède par récurrence sur $n$. 

Pour $n=p+1$, il suffit de prendre $\varphi=0$ et de choisir pour $g$
la métrique donnée par le lemme~\ref{multi:lem1}. En effet, pour
$\varphi=0$, on a $\tilde\mu_{q,i}(S^n,g,\varphi)=\mu_{q,i}(S^n,g)
=\mu_{n-q-1,i}(S^n,g)$ pour tout $q$. En particulier, on a
$\tilde\mu_{p,i}(S^n,g,\varphi)=\mu_{n-1,i}(S^n,g)=\mu_{0,i}(S^n,g)$, et pour
les autres degrés le spectre est minoré par la même constante $C$.

Considérons maintenant un entier $n>p+1$ et supposons le lemme vrai
pour la dimension $n-1$. On a donc sur $S^{n-1}$ une métrique $\bar g$
et une fonction $\bar\varphi$ satisfaisant la conclusion du lemme. 
On considère un domaine $U\subset S^n$ qui est le produit de 
$(S^{n-1},\bar g,\bar\varphi)$ avec un petit intervalle $[-\varepsilon,
\varepsilon]$, on prolonge la métrique et la fonction de manière quelconque
sur $S^n$. Pour $\varepsilon$ suffisamment petit, le début du spectre
sur $U$ est le même que celui de $(S^{n-1},\bar g,\bar\varphi)$, avec
la même propriété de stabilité d'après le lemme~\ref{multi:produit}.
On applique alors le théorème~\ref{conv:th2} avec $\alpha=\frac n2-1$,
de manière à avoir convergence du spectre pour les degrés $0\leq q<n-1$
(comme la cohomologie du domaine $U$ est triviale sauf en degré $0$ et $n-1$
et que $H^0(S^n/U)$ est aussi trivial, il n'y a pas de valeurs propres qui 
tendent vers 0 pour ces degrés).
L'argument de stabilité fournit alors une métrique $g$ et une
fonction $\varphi$ telles que $\tilde\mu_{p,i}(S^n,g,\varphi)=a_i$ pour
$i\leq N$ (avec stabilité), $\tilde\mu_{p,N+1}(S^n,g,\varphi)>C$ et 
$\tilde\mu_{q,1}(S^n,g,\varphi)>C$ pour $1\leq q< n-1$, $q\neq p$. Comme
le volume de $U$ peut être choisit arbitrairement petit, le résultat
en degré $q=n-1$ découle du point 3 du théorème~\ref{conv:th2}.
\end{demo}
Comme dans \cite{ja11}, le cas des 1-formes est traité séparément:
\begin{lem}\label{multi:lem3}
Pour tous entiers $N\geq1$, $n\geq4$, toute suite finie
$0<a_1\leq a_2\leq\ldots\leq a_N$ et toute constante $C>a_N$, il existe
une métrique $g$ et une fonction $\varphi$ sur la variété
$S^n$ telles que $\tilde\mu_{1,i}(S^n,g,\varphi)=a_i$ pour
$i\leq N$, ces valeurs propres vérifiant l'hypothèse de stabilité,
$\tilde\mu_{1,N+1}(S^n,g,\varphi)>C$ et $\tilde\mu_{p,1}(S^n,g,\varphi)>C$ 
pour $p>1$, le volume $\Vol(S^n,g)$ étant arbitrairement petit.
\end{lem}
\begin{demo}
On applique le lemme précédent pour $p=n-2$. Il suffit ensuite de changer
le signe de $\varphi$, la relation (\ref{theo:hodge}) permet de conclure.
\end{demo}

On déduit des lemmes précédents un modèle de multiplicité pour
chaque degré $p$ sur des boules de dimension $n$:

\begin{lem}\label{multi:lem4}
Pour tous entiers $N\geq1$, $n\geq4$ et $p\in[1,n-1]$, toute suite finie
$0<a_1\leq a_2\leq\ldots\leq a_N$ et toute constante $C>a_N$, il existe
une métrique $g$ et une fonction $\varphi$ sur la boule
$B^n$ telles que $\tilde\mu_{p,i}(B^n,g,\varphi)=a_i$ pour
$i\leq N$, ces valeurs propres vérifiant l'hypothèse de stabilité,
$\tilde\mu_{p,N+1}(B^n,g,\varphi)>C$ et $\tilde\mu_{q,1}(B^n,g,\varphi)>C$
pour $q\neq p$, le volume $\Vol(S^n,g)$ étant arbitrairement petit.
\end{lem}
\begin{demo}
Il suffit d'appliquer le théorème~\ref{boule:theo} et l'argument de 
stabilité aux modèles construits précédemment: on utilise le 
lemme~\ref{multi:lem2} pour les degrés $2\leq p\leq n-1$ et le 
lemme~\ref{multi:lem3} pour le degré 1.
\end{demo}

On va maintenant montrer les théorèmes~\ref{intro:th1} et~\ref{intro:th2}

\begin{demo}[du théorème~\ref{intro:th1}]
On commence par se donner, pour chaque degré $p$, un domaine $U_p$ de $M$
(muni d'une métrique et d'une fonction lisse)
dont le début du spectre est celui qu'on veut prescrire. On applique 
pour cela le lemme~\ref{multi:lem4} avec $C>\sup_{p,i}a_{p,i}$. On
prolonge la métrique et la fonction à $M$.

On applique ensuite le théorème~\ref{conv:th2} au domaine $U=\bigcup U_p$
avec $\alpha>\frac n2$. Comme $U$ n'a de cohomologie qu'en degré~0, 
on a convergence du spectre pour les autres degrés. L'argument de stabilité
assure alors l'existence d'une métrique $g$ et d'une fonction $\varphi$
sur $M$ telle que $\tilde\mu_{p,i}(M,g,\varphi)=a_{p,i}$ pour tout
$p\geq1$ et tout $1\leq i\leq N$

Selon un argument classique (utilisé dans \cite{gu04}, \cite{ja08} et
\cite{ja11}), on peut prescrire simultanément le volume
en le traitant comme une valeur propre simple : on ajoute
à la famille $U_p$ une boule de volume $V-\sum_p\Vol(U_p)$
(en ayant choisit les volumes des $U_p$ suffisamment petits) et dont
les valeurs propres pour les $p$-formes sont arbitrairement grandes
(c'est possible selon \cite{gp95}). Le fait que le volume de cette boule
vérifie l'hypothèse de transversalité signifie simplement qu'on
peut lui donner n'importe quelle valeur au voisinage de $V-\sum_p\Vol(U_p)$,
par exemple par homothétie.
\end{demo}

\begin{demo}[du théorème~\ref{intro:th2}]
Soit $\Sigma$ une surface compacte. Y.~Colin de Verdière a montré dans 
\cite{cdv86} (corollaire~7.2) que pour toute métrique $g$ sur $\Sigma$, 
il existe une fonction $\varphi$ telle que la multiplicité de
$\tilde\mu_{0,1}(\Sigma,g,\varphi)$ soit de multiplicité stable 
$C(\Sigma)-1$. En appliquant la formule (\ref{theo:hodge}), on
obtient que $\tilde\mu_{1,1}(\Sigma,g,-\varphi)$ a la même multiplicité.

On peut plonger la surface $\Sigma$ dans $M$ et considérer le voisinage
$U=\Sigma\times[-\varepsilon,\varepsilon]$ muni de la métrique produit
qu'on notera encore $g$ et en relevant la fonction $\varphi$ à $U$.
Pour $\varepsilon$ suffisamment petit, $\tilde\mu_{1,1}(U,g,-\varphi)$
est encore de multiplicité stable $C(\Sigma)-1$ d'après le 
lemme~\ref{multi:produit}. 

 Après avoir étendu $g$ et $\varphi$ à $M$ de manière quelconque, on peut 
appliquer le théorème de convergence~\ref{conv:th2} à $U$ avec $\alpha>2$.
Si $H^1(M)\to H^1(\Sigma)$ est surjective, on a alors $H^1(M/U)=\{0\}$ et 
l'argument de stabilité donne le résultat.
\end{demo}

\section*{Appendice}
\renewcommand{\theequation}{A.\arabic{equation}}

Le but de cet appendice est de justifier l'équivalence des trois formulations
suivantes du laplacien de Witten :
\begin{theo}\label{app:th1}
Pour toute fonction $\varphi$, on a
\begin{eqnarray*}
\tilde\Delta_\varphi&=&\Delta+|\de\varphi|^2+\mathcal L_{\nabla\varphi}+
\mathcal L_{\nabla\varphi}^*\\
&=&\Delta+|\de\varphi|^2+\Delta\varphi-(\mathcal L_Xg_p)\\
&=&\Delta+|\de\varphi|^2+\Delta\varphi+2(\Hess\varphi)
\end{eqnarray*}
\end{theo}
On va en fait montrer un résultat plus général. Si $X$ est un champ de 
vecteur, on pose 
\begin{equation}\label{app:eq1}
\tilde\de_X\omega=\de\omega+X^\flat\wedge\omega\textrm{ et }
\tilde\codiff_X\omega=\codiff\omega+\prodint_X\omega.
\end{equation}
On retrouve les définitions de $\tilde\de_\varphi$ et $\tilde\codiff_\varphi$
en remplaçant $X$ par $\nabla\varphi$. Le théorème~\ref{app:th1} est
un cas particulier du résultat qui suit. Si $\alpha$ est une $1$-forme, 
$\overline\nabla\alpha$ la partie symétrique de la dérivée covariante
$\nabla\alpha$ qu'on identifie à un endomorphisme de $\Lambda^1TM$ :
\begin{theo}\label{app:th2}
Pour tout champ de vecteur $X$, on a
\begin{eqnarray*}
\tilde\codiff_X\tilde\de_X+\tilde\de_X\tilde\codiff_X
&=&\Delta+|X|^2+\mathcal L_X+\mathcal L_X^*\\
&=&\Delta+|X|^2+\mathrm{div} X-(\mathcal L_Xg_p)\\
&=&\Delta+|X|^2+\mathrm{div} X+2(\overline\nabla X^\flat)
\end{eqnarray*}
\end{theo}
Comme dans le théorème~\ref{app:th1}, l'action de $\overline\nabla X^\flat$
s'étend naturellement aux $p$-formes (cf. équation~(\ref{theo:eq4})),
dans la suite on notera $(\overline\nabla X^\flat)_p$ cette extension
s'il est nécessaire de préciser le degré.
Elle s'exprime en fait de manière très simple si on utilise le 
formalisme de l'algèbre de courbure. Comme ce formalisme va intervenir
dans la démonstration, nous allons le rappeler.

On note $C^p=S^2\Lambda^pTM$ l'espace des formes quadratiques sur les
$p$-formes. Ses éléments peuvent aussi s'interpréter comme endomorphismes 
de $\Lambda^pTM$ ou $2p$-tenseurs antisymétriques par rapport aux $p$ 
premières variables et symétriques par échanges des $p$ premières et des
$p$ dernières. L'espace des tenseurs de courbure est défini par $C^*=
\oplus_{p\geq0}C^p$. Il contient entre autres les tenseurs de riemann et 
de Weyl ($p=2$) la courbure de Ricci ($p=1$) et la courbure scalaire,
ce qui justifie son nom. Il est muni d'une structure d'algèbre commutative 
dont le produit $\owedge$ est défini par $\omega\otimes\omega\owedge
\theta\otimes\theta=\omega\wedge\theta\otimes\omega\wedge\theta$. Avec
ces notations, on a $g_p=\frac1{p!}g_1^p$ et 
$(\overline\nabla X^\flat)_p=(\overline\nabla X^\flat)_1\owedge g_{p-1}$
(voir par exemple \cite{la05}, \cite{lahdr} et les références qui y sont
données).
\begin{remarque}
Pour un champ de vecteur $X$ quelconque, on n'a pas $\tilde\de_X^2=0$,
et par conséquent $(\tilde\de_X+\tilde\codiff_X)^2$ est différent de
$\tilde\codiff_X\tilde\de_X+\tilde\de_X\tilde\codiff_X$ en général.
\end{remarque}

La première égalité du théorème~\ref{app:th2} se calcule de manière 
immédiate, comme celle du théorème~\ref{app:th1}: à partir de 
(\ref{app:eq1}), on obtient
\begin{equation}
\tilde\codiff_X\tilde\de_X+\tilde\de_X\tilde\codiff_X=\Delta+|X|^2+
\prodint_X\de+\de\prodint_X+X^\flat\wedge\codiff+\codiff(X^\flat\wedge).
\end{equation}
Or, $\prodint_X\de+\de\prodint_X=\mathcal L_X$ (formule de Cartan), et
$\prodint_X\de+\de\prodint_X$ est l'adjoint de $\codiff(X^\flat\wedge)+
X^\flat\wedge\codiff$. Les autres égalités découlent du 
\begin{lem}
Pour tout un champ de vecteur $X$ et toute $p$-forme $\omega$ sur $M$,
on a l'identité
$$\mathcal L_X\omega+\mathcal L_X^*\omega
=(\mathrm{div} X)\omega-(\mathcal L_Xg_p)\omega
=(\mathrm{div} X)\omega+2(\overline\nabla X^\flat)\omega$$.
\end{lem}

\begin{demo}
Dans toute la démonstration on supposera que la variété est compacte
sans bord, pour simplifier les intégrations par parties. Le cas à bord
s'en déduit par restriction à un domaine.

Soit $\omega$ et $\omega'$ deux $p$-formes sur $M$. En dérivant selon $X$ 
leur produit scalaire, il vient
\begin{equation}
X\cdot\langle\omega,\omega'\rangle=(\mathcal L_Xg_p)(\omega,\omega')+
\langle\mathcal L_X\omega,\omega'\rangle
+\langle\omega,\mathcal L_X\omega'\rangle.
\end{equation}
On intègre les différents termes de cette équation. Le membre de 
gauche donne
\begin{equation}
\int_M X\cdot\langle\omega,\omega'\rangle=\int_M\langle X,\nabla
\langle\omega,\omega'\rangle\rangle
=\int_M(\mathrm{div} X)\langle\omega,\omega'\rangle
\end{equation}
et le dernier terme devient
\begin{equation}
\int_M\langle\omega,\mathcal L_X\omega'\rangle=\int_M\langle\mathcal L_X^*
\omega,\omega'\rangle.
\end{equation}
On obtient
\begin{equation}
\int_M(\mathrm{div} X)\langle\omega,\omega'\rangle=
(\mathcal L_Xg_p)(\omega,\omega')+\int_M\langle\mathcal L_X\omega,\omega'\rangle
+\langle\mathcal L_X^*\omega,\omega'\rangle.
\end{equation}
Comme cette égalité est vraie pour tout $\omega'$, on en déduit l'égalité
ponctuelle
\begin{equation}
\mathcal L_X\omega+\mathcal L_X^*\omega=
(\mathrm{div} X)\omega-(\mathcal L_Xg_p)\omega.
\end{equation}

On doit réexprimer $\mathcal L_Xg_p$ pour obtenir la deuxième expression
du lemme. On va commencer par le cas $p=1$. On note $U$ et $V$ deux champs 
de vecteurs quelconques et $\alpha$, $\beta$ leur 1-forme duale. En utilisant 
le fait que $\langle U,V\rangle=\langle\alpha,\beta\rangle=\alpha(V)=\beta(U)$
et les relations de dérivations suivantes :
\begin{eqnarray}\label{app:eq2}
X\cdot\langle\alpha,\beta\rangle & = & (\mathcal L_Xg_1)(\alpha,\beta)+
\langle\mathcal L_X\alpha,\beta\rangle+
\langle\alpha,\mathcal L_X\beta\rangle\nonumber\\
&=&(\mathcal L_Xg_1)(\alpha,\beta)+\mathcal L_X\alpha(V)+\mathcal L_X\beta(U)
\end{eqnarray}
\begin{equation}\label{app:eq3}
X\cdot\alpha(V)=\mathcal L_X\alpha(V)+\alpha([X,V])
\end{equation}
\begin{equation}\label{app:eq4}
X\cdot\beta(U)=\mathcal L_X\beta(U)+\beta([X,U])
\end{equation}
et en utilisant le fait que la dérivée covariante est sans torsion:
\begin{eqnarray}\label{app:eq5}
\nabla_X\langle U,V\rangle & = &
\langle\nabla_XU,V\rangle+\langle U,\nabla_XV\rangle\nonumber\\
& = &\langle\nabla_UX,V\rangle+\langle[X,U],V\rangle
+\langle U,\nabla_VX\rangle+\langle U,[X,V]\rangle\nonumber\\
& = & \langle\nabla_UX,V\rangle+\beta([X,U])+\langle U,\nabla_VX\rangle
+\alpha([X,U]),
\end{eqnarray}
on obtient, en simplifiant la combinaison $\mathrm{(\ref{app:eq2})}+
\mathrm{(\ref{app:eq5})}-\mathrm{(\ref{app:eq3})}-\mathrm{(\ref{app:eq4})}$:
\begin{equation}
(\mathcal L_Xg_1)(\alpha,\beta)+\langle\nabla_UX,V\rangle+
\langle U,\nabla_VX\rangle=0,
\end{equation}
soit
\begin{equation}
(\mathcal L_Xg_1)(\alpha,\beta)=-2(\overline\nabla X)(U,V).
\end{equation}
Par dualité, on obtient bien que $(\mathcal L_Xg_1)\omega=
-2(\overline\nabla X^\flat)\omega$ si $\omega$ est une 1-forme.

Pour étendre se résultat aux $p$-formes il suffit d'utiliser le fait
que $g_p=\frac1{p!}g_1^p$ et d'appliquer la formule de Leibniz
\begin{equation}
(\mathcal L_Xg_p)=(\mathcal L_X\frac1{p!}g_1^p)=
(\mathcal L_Xg_1)\owedge\frac1{(p-1)!}g_1^{p-1}=
-2(\overline\nabla X^\flat)\owedge g_{p-1},
\end{equation}
Ce qui achève la démonstration.
\end{demo}

\noindent Pierre \textsc{Jammes}\\
Université d'Avignon et des pays de Vaucluse\\
Laboratoire d'analyse non linéaire et géométrie (EA 2151)\\
F-84018 Avignon\\
\texttt{Pierre.Jammes@ens-lyon.org}

\begin{thebibliography}{CdVT93}
{\small
\makeatletter
\ifx\fonteauteurs\@undefined
\def\fonteauteurs{\scshape}\fi
\makeatother

\bibitem[AC93]{ac93}
\bgroup\fonteauteurs C.~Ann\'e\egroup{} et \bgroup\fonteauteurs
  B.~Colbois\egroup{} -- \og~Op\'erateur de hodge-laplace sur des vari\'et\'es
  compactes priv\'ees d'un nombre fini de boules~\fg, {\em J. Funct. Anal.},
  115, p.~190--211, 1993.

\bibitem[An89]{an89}
\bgroup\fonteauteurs C.~Anné\egroup{} -- \og~Principe de {D}irichlet pour les
  formes différentielles~\fg, {\em Bull. soc. math. France}, 117 (4),
  p.~445--450, 1989.

\bibitem[BCC98]{bcc98}
\bgroup\fonteauteurs G.~Besson\egroup{}, \bgroup\fonteauteurs
  B.~Colbois\egroup{} et \bgroup\fonteauteurs G.~Courtois\egroup{} -- \og~Sur
  la multiplicité de la première valeur propre de l'opérateur de {S}chrödinger
  avec champ magnétique sur la sphère ${S}^2$~\fg, {\em Trans. Amer. Math.
  Soc.}, 350 (1), p.~331--345, 1998.

\bibitem[Be80]{be80}
\bgroup\fonteauteurs G.~Besson\egroup{} -- \og~Sur la multiplicité de la
  première valeur propre des surfaces riemanniennes~\fg, {\em Ann. inst.
  Fourier}, 30 (1), p.~109--128, 1980.

\bibitem[CdV86]{cdv86}
\bgroup\fonteauteurs Y.~Colin~de Verdière\egroup{} -- \og~Sur la multiplicité
  de la première valeur propre non nulle du laplacien~\fg, {\em Comment. Math.
  Helv.}, 61 (2), p.~254--270, 1986.

\bibitem[CdV87]{cdv87}
\bgroup\fonteauteurs Y.~Colin~de Verdière\egroup{} -- \og~Construction de
  laplaciens dont une partie finie du spectre est donnée~\fg, {\em Ann. scient.
  \'Ec. norm. sup.}, 20 (4), p.~99--615, 1987.

\bibitem[CdV88]{cdv88}
\bgroup\fonteauteurs Y.~Colin~de Verdière\egroup{} -- \og~Sur une hypothèse de
  transversalité d'{A}rnol'd~\fg, {\em Comment. Math. Helv.}, 63 (2),
  p.~184--193, 1988.

\bibitem[CdVT93]{cdvt93}
\bgroup\fonteauteurs Y.~Colin~de Verdière\egroup{} et \bgroup\fonteauteurs
  N.~Torki\egroup{} -- \og~Opérateur de {S}chrödinger avec champ
  magnétique~\fg, {\em Sémin. théor. spectr. géom.}, 11, p.~9--18, 1993.

\bibitem[Ch76]{ch76}
\bgroup\fonteauteurs S.~Y. Cheng\egroup{} -- \og~Eigenfunctions and nodal
  sets~\fg, {\em Comment. Math. Helv.}, 51 (1), p.~43--55, 1976.

\bibitem[Da05]{da05}
\bgroup\fonteauteurs M.~Dahl\egroup{} -- \og~Prescribing eigenvalues of the
  {D}irac operator~\fg, {\em Manuscripta Math.}, 118 (2), p.~191--199, 2005,
  math.DG/0311172.

\bibitem[Do82]{do82}
\bgroup\fonteauteurs J.~Dodziuk\egroup{} -- \og~Eigenvalues of the {L}aplacian
  on forms~\fg, {\em Proc. of Am. Math. Soc.}, 85, p.~438--443, 1982.

\bibitem[Er02]{er02}
\bgroup\fonteauteurs L~Erd{\H{o}}s\egroup{} -- \og~Spectral shift and
  multiplicity of the first eigenvalue of the magnetic {S}chr\"odinger operator
  in two dimensions~\fg, {\em Ann. inst. Fourier}, 52 (6), p.~1833--1874, 2002.

\bibitem[GP95]{gp95}
\bgroup\fonteauteurs G.~Gentile\egroup{} et \bgroup\fonteauteurs
  V.~Pagliara\egroup{} -- \og~Riemannian metrics with large first eigenvalue on
  forms of degree $p$~\fg, {\em Proc. of Am. Math. Soc.}, 123 (12),
  p.~3855--3858, 1995.

\bibitem[GT06]{gt06}
\bgroup\fonteauteurs V.~Gold'shtein\egroup{} et \bgroup\fonteauteurs
  M.~Troyanov\egroup{} -- \og~Sobolev inequalities for differential forms ans
  ${L}_{q,p}$-cohomology~\fg, {\em J. Geom. Anal.}, 16 (4), p.~597--632, 2006,
  math.DG/0506065.

\bibitem[Gu04]{gu04}
\bgroup\fonteauteurs P.~Guérini\egroup{} -- \og~Prescription du spectre du
  laplacien de {H}odge-de~{R}ham~\fg, {\em Ann. scient. \'Ec. norm. sup.}, 37
  (2), p.~270--303, 2004.

\bibitem[He85]{he84}
\bgroup\fonteauteurs G.~Henniart\egroup{} -- \og~Les in\'egalit\'es de
  {M}orse~\fg, Dans {\em {S}\'eminaire Bourbaki 83/84},  volume 121--122 de
  {\em Ast\'erisque},  pages  43--61, Soc. Math. France, 1985.

\bibitem[HN06]{hn06}
\bgroup\fonteauteurs B.~Heffler\egroup{} et \bgroup\fonteauteurs
  F.~Nier\egroup{} -- \og~Quantitative analysis of metastability in reversible
  diffusion processes via a {W}itten complex approach: the case with
  boundary~\fg, {\em Mem. soc. math. France}, 105, 2006.

\bibitem[HHN99]{hohon99}
\bgroup\fonteauteurs M.~Hoffmann-Ostenhof\egroup{}, \bgroup\fonteauteurs
  T.~Hoffmann-Ostenhof\egroup{} et \bgroup\fonteauteurs
  N.~Nadirashvili\egroup{} -- \og~On the multiplicity of eigenvalues of the
  {L}aplacian on surfaces~\fg, {\em Ann. Global Anal. Geom.}, 17 (1),
  p.~43--48, 1999.

\bibitem[Ja07a]{ja07a}
\bgroup\fonteauteurs P.~Jammes\egroup{} -- \og~Extrema de valeurs propres dans
  une classe conforme~\fg, {\em S\'emin. th\'eor. spectr. g\'eom.}, 24,
  p.~23--43, 2007, arXiv:0804.0488.

\bibitem[Ja07b]{ja07b}
\bgroup\fonteauteurs P.~Jammes\egroup{} -- \og~Minoration conforme du spectre
  du laplacien de {H}odge-de~{R}ham~\fg, {\em Manuscripta Math.}, 123 (1),
  p.~15--23, 2007, math.DG/0604591.

\bibitem[Ja08]{ja08}
\bgroup\fonteauteurs P.~Jammes\egroup{} -- \og~Prescription du spectre du
  laplacien de {H}odge-de~{R}ham dans une classe conforme~\fg, {\em Comment.
  Math. Helv.}, 83 (3), p.~521--537, 2008, math.DG/0601738.

\bibitem[Ja09]{ja09}
\bgroup\fonteauteurs P.~Jammes\egroup{} -- \og~Construction de valeurs propres
  doubles du laplacien de {H}odge-de~{R}ham~\fg, {\em J. Geom. Anal.}, 19 (3),
  p.~643--654, 2009, math.DG/0608758.

\bibitem[Ja11]{ja11}
\bgroup\fonteauteurs P.~Jammes\egroup{} -- \og~Prescription de la
  multiplicit\'e des valeurs propres du laplacien de {H}odge-de~{R}ham~\fg,
  {\em Comment. Math. Helv.}, 86 (4), p.~967--984, 2011, arXiv:0804.0104.

\bibitem[La05]{la05}
\bgroup\fonteauteurs M.-L. Labbi\egroup{} -- \og~{D}ouble forms, curvature
  structures and the $(p,q)$ curvatures~\fg, {\em Trans. Amer. Math. Soc.}, 357
  (10), p.~3971--3992, 2005, math.DG/0404081.

\bibitem[La06]{lahdr}
\bgroup\fonteauteurs M.-L. Labbi\egroup{} -- {\em Courbure riemannienne :
  diff\'erentes notions de positivit\'e}, habilitation \`a diriger des
  recherches, universit\'e de Montpellier {II}, 2006, math.DG/0611371.

\bibitem[Lo96]{lo96}
\bgroup\fonteauteurs J.~Lohkamp\egroup{} -- \og~Discontinuity of geometric
  expansions~\fg, {\em Comment. Math. Helv.}, 71 (2), p.~213--228, 1996.

\bibitem[{Mc}93]{mc93}
\bgroup\fonteauteurs J.~{Mc}{Gowan}\egroup{} -- \og~The $p$-spectrum of the
  {L}aplacian on compact hyperbolic three manifolds~\fg, {\em Math. Ann.}, 297
  (4), p.~725--745, 1993.

\bibitem[Na88]{na88}
\bgroup\fonteauteurs N.~Nadirashvili\egroup{} -- \og~Multiple eigenvalues of
  the {L}aplace operator~\fg, {\em Math. USSR-Sb.}, 61 (1), p.~225--238, 1988.

\bibitem[Sé94]{se94}
\bgroup\fonteauteurs B.~Sévennec\egroup{} -- \og~Multiplicité du spectre des
  surfaces: une approche topologique~\fg, {\em Sémin. théor. spectr. géom.},
  12, p.~29--36, 1994.

\bibitem[Sé02]{se02}
\bgroup\fonteauteurs B.~Sévennec\egroup{} -- \og~Multiplicity of the second
  {S}chrödinger eigenvalue on closed surfaces~\fg, {\em Math. Ann.}, 324 (1),
  p.~195--211, 2002.

\bibitem[Ta96]{ta96}
\bgroup\fonteauteurs M.~Taylor\egroup{} -- {\em Partial differential equations
  I}, Springer, 1996.

\bibitem[Wi82]{wi82}
\bgroup\fonteauteurs E.~Witten\egroup{} -- \og~Supersymmetry and {M}orse
  theory~\fg, {\em J. Differ. Geom.}, 17, p.~661--692, 1982.

}\end{thebibliography}
\end{document}